\documentclass[12pt]{article}
\usepackage[T2A]{fontenc}
\usepackage[cp866]{inputenc}
\usepackage[english,russian]{babel}
\usepackage[tbtags]{amsmath}
\usepackage{amsfonts,amssymb,mathrsfs,amscd,amsthm,hyper}
\usepackage{mathtools}
\RequirePackage{graphicx}
\theoremstyle{plain}

\theoremstyle{plain}

\newtheorem{theorem}{Теорема}
\newtheorem{theoremA}{Теорема~A}
\newtheorem{theoremSt}{Теорема Шталя}
\newtheorem{theoremFsSt}{Первая теорема Шталя}
\newtheorem{theoremSeSt}{Вторая теорема Шталя}
\newtheorem{theoremJS}{Теорема Йенча--Сегё}
\newtheorem{lemma}{Лемма}

\theoremstyle{definition}
\newtheorem{definition}{Определение}
\newtheorem{remark}{Замечание}
\let\ge\geqslant\let\leq\leqslant\let\geq\geqslant
\def\supp{\operatorname{supp}}
\def\mdeg{\operatorname{deg}}
\let\eps\varepsilon
\let\myo\overline\let\myh\widehat\let\pfi\varphi\let\myt\widetilde
\let\leq\leqslant
\let\geq\geqslant
\def\({\left(}
\def\){\right)}
\def\mM{M}
\def\EE{E}
\def\mymu{\mu}

\def\Re{\operatorname{Re}}
\def\Im{\operatorname{Im}}
\def\mdeg{\operatorname{deg}}

\def\hh{h}
\def\const{\operatorname{const}}

\def\mcap{\operatorname{cap}}

\def\Zh{\operatorname{Zh}}
\def\Res{\operatornamewithlimits{Res}}
\def\DD{\mathbb D}
\def\EE{E}
\def\mM{M}

\def\zz{\mathbf z}

\def\sU{U}
\def\myA{\mathscr A}
\def\HH{\mathscr H}
\def\KK{\mathscr K}
\def\sV{\mathscr V}
\def\sR{\mathscr R}
\def\SR{\mathscr R}
\def\RS{\mathfrak R}
\def\sH{\mathscr H}
\def\sM{\mathscr M}
\def\MM{\mathscr M}
\def\FF{\mathscr F}
\def\sG{\mathscr G}
\def\RR{\mathbb R}
\def\CC{\mathbb C}
\def\NN{\mathbb N}
\def\RS{\mathfrak R}
\def\rH{\mathrm H}

\def\bad{\spaceskip=0.33emplus0.6emminus0.15em\immediate\write5{\string\bad}}
\begin{document}

\title{Сходимость нелинейных аппроксимаций Паде--Чебышёва для многозначных
аналитических функций, вариация равновесной энергии и $S$-свойство
стационарных компактов}

\author{А.\,А.~Гончар, Е.\,А.~Рахманов, С.\,П.~Суетин}

\maketitle

\begin{abstract}
{\bad
Исследуется сходимость диагональных нелинейных аппроксимаций Паде--Чебышёва для
многозначных аналитических функций, заданных и вещественнозначных на единичном
отрезке $[-1,1]$. Получен аналог классической теоремы Шталя с сходимости по
емкости таких аппроксимаций в соответствующей ``максимальной'' области
мероморфности заданной функции. Скорость сходимости этих рациональных функций
на отрезке и в максимальной области мероморфности заданной
функции характеризуется в терминах ``стационарного'' компакта для
некоторой смешанной (гриново-логарифмической) теоретико-потенциальной задачи
равновесия. Найдено предельное распределение точек интерполяции на отрезке
$[-1,1]$ исходной функции нелинейной аппроксимацией Паде--Чебышёва.
}

Библиография: 57 названий.
\end{abstract}

\markright{Аппроксимации Паде--Чебышёва}

\footnotetext[0]{Работа выполнена при поддержке Российского фонда
фундаментальных исследований (гранты \No~08-01-00317 и
\No~09-01-12160-офи-м) и Программы поддержки ведущих научных школ РФ (грант
\No~НШ-8033.2010.1).}

\section{Введение}\label{s1}

\subsection{}\label{ss1}

Пусть на единичном отрезке $E=[-1,1]$ задана вещественная функция~$f$, голоморфная
на~$E$. Пусть $\{T_k(x)\}$ -- полиномы Чебышёва (первого рода), ортонормированные на
отрезке $E$ по мере $d\tau(x)=(1-x^2)^{-1/2}\,dx$. Функции~$f$
соответствует ряд Фурье--Чебышёва:
\begin{equation}
f(x)\sim\sum_{k=0}^\infty c_kT_k(x),\qquad
c_k=c_k(f)=\int_{-1}^1 f(x)T_k(x)\,d\tau(x).
\label{1.01}
\end{equation}
Для произвольного $\rho>1$ через $\Gamma_\rho$ обозначим эллипс с фокусами
в точках $\pm1$ и суммой полуосей, равной $\rho$. Внутренность эллипса
$\Gamma_\rho$ обозначим через $D_\rho$; области $D_\rho$ будем называть
{\it каноническими} (относительно $E$).
Так как функция $f$ голоморфна на $E$ ($f\in\HH(E)$), то ряд~\eqref{1.01}
сходится к $f(x)$ локально равномерно в области $D(f)=D_\rho$, где
$\rho=\rho_0(f)$ -- индекс максимальной канонической области, в которую $f$
продолжается как голоморфная функция. Область $D(f)$ --
максимальная область сходимости ряда~\eqref{1.01}.

Для произвольного натурального числа $n$ через $\SR_n$ обозначим класс всех
рациональных функций $r$ вида $r=p/q$, где $p,q$ -- произвольные алгебраические
полиномы степени $\leq{n}$, $q\not\equiv0$. Пусть функция $f$ продолжается с отрезка
$E$ в область $G\supset{E}$ как (однозначная) мероморфная функция. Тогда
для величины
\begin{equation}
\eps_n(f)=\inf_{r\in\SR_n}\|f-r\|_{E}
\label{1.02}
\end{equation}
-- наилучшего равномерного приближения $f$ рациональными дробями класса $\SR_n$
-- справедливо {\it неравенство Уолша} (см.~\cite{Wal61}):
\begin{equation}
\varlimsup_{n\to\infty}\eps_n^{1/n}\leq{q},
\label{1.03}
\end{equation}
где
\begin{equation}
q=\exp\biggl\{-\frac1{C(E,\myo{\CC}\setminus{G})}\biggr\}<1,
\label{1.031}
\end{equation}
$C(E,\myo{\CC}\setminus{G})$ -- емкость конденсатора
$(E,\myo{\CC}\setminus{G})$.

Неравенство~\eqref{1.03} носит {\it универсальный
характер} -- оно неулучшаемо в классе {\it всех} функций, мероморфных в
области $G$. Однако для некоторых вполне естественных классов функций
неравенство Уолша~\eqref{1.03} может быть существенно усилено.

\normalbaselineskip=12.pt %

Пусть $f=\myh\sigma$, где $\myh\sigma$ -- {\it марковская функция}:
\begin{equation}
\myh{\sigma}(z)=\int_c^d\frac{d\sigma(x)}{z-x},
\quad z\in\myo{\CC}\setminus{[c,d]},
\label{1.04}
\end{equation}
$\sigma$ -- положительная борелевская мера на отрезке $[c,d]\subset\RR$,
$[c,d]\cap{E}=\varnothing$. Если $\sigma'=d\sigma/dx>0$ почти всюду (п.в.) на
$[c,d]$, то для наилучших рациональных аппроксимаций функции $f$ имеем
(см.~\cite{Gon78})
\begin{equation}
\lim_{n\to\infty}\eps_n^{1/n}={q^2},
\label{1.05}
\end{equation}
где величина $q<1$ имеет тот же смысл, что и в~\eqref{1.03} с
заменой конденсатора $(E,\myo{\CC}\setminus{G})$ на конденсатор $(E,[c,d])$
(в соотношении~\eqref{1.05} содержатся два утверждения: существование предела
и информация о его величине).
Таким образом, в классе марковских функций~\eqref{1.04}, голоморфных в области
$G=\myo{\CC}\setminus[c,d]$, вместо неравенства Уолша~\eqref{1.03} имеет
место существенно более сильное {\it равенство}~\eqref{1.05}.
Соотношение~\eqref{1.05} справедливо и для функций вида $f=\myh\sigma+r$, где
$r$ -- рациональная функций, голоморфная на $E$.

В серии работ 1985--1986~гг. Г.~Шталь~\cite{Sta85a}--\cite{Sta86b} (см.
также~\cite{Sta95},~\cite{Sta97} и
приложение~\ref{apl}) получил ряд результатов о
сходимости классических и многоточечных аппроксимаций Паде для функций,
голоморфных соответственно в точке $z=\infty$ и на фиксированном односвязном
континууме в $\CC$ и имеющих конечное число особых точек многозначного характера.
Из его результата, относящегося к случаю континуума, и результатов
работы~\cite{GoRa87} вытекает справедливость {\it равенства}~\eqref{1.05} в классе
функций, вещественных голоморфных на~$E$ и имеющих вне $E$ конечное число
особых точек многозначного характера.
Точнее, обозначим через $R_n$ наилучшие равномерные рациональные аппроксимации
функции $f$ на отрезке $E$ в классе $\SR_n$. Тогда справедлива следующая

\begin{theoremSt}[(см.~\cite{Sta86a}--\cite{Sta86b} и приложение~\ref{apl})]
Пусть функция $f$ вещественнозначна и голоморфна на отрезке $E=[-1,1]$ и допускает
мероморфное продолжение с отрезка в $\myo{\CC}$ всюду за исключением конечного
числа особых точек многозначного характера.
Тогда существует единственный компакт $F=F(f)$, не пересекающийся с $E$ и
такой, что функция $f$ мероморфна\footnote{Точнее, функция $f$ допускает
мероморфное продолжение в $\myo{\CC}\setminus{F}$, которое мы также будем
обозначать через~$f$.}
в области $\myo{\CC}\setminus{F}$ и для любого компакта
$K\subset\myo{\CC}\setminus(\EE\cup{F})$
\begin{equation}
\bigl|(f-R_n)(z)\bigr|^{1/n}\overset{\mcap}{\longrightarrow}
e^{-2G_{F}^{\lambda}(z)}<1,
\qquad z\in K.
\label{1.06}
\end{equation}
\end{theoremSt}
В соотношении~\eqref{1.06} $\lambda$ -- единичная мера с носителем на $E$,
$G_{F}^\lambda(z)$ -- гринов (относительно $F$) потенциал меры $\lambda$.
Компакт $F$ обладает так называемым $S$-свойством (см.~\eqref{spro1}), не
разбивает плоскость и состоит из конечного числа кусочно-аналитических дуг, а
$\lambda$ -- равновесная мера (подробнее см.~\S\,\ref{s2}).
Непосредственно из~\eqref{1.06} вытекает, что функция $f$ продолжается с
отрезка $E$ в $\myo{\CC}\setminus{F}$ как однозначная мероморфная функция.
С учетом результатов работы~\cite{GoRa87} из~\eqref{1.06} вытекает, что
в условиях и обозначениях теоремы Шталя справедлив аналог
равенства~\eqref{1.05}
\begin{equation}
\lim_{n\to\infty}\eps_n^{1/n}={q^2};
\label{1.07}
\end{equation}
величина $q$ имеет тот же смысл, что и в~\eqref{1.03} с
заменой конденсатора $(E,\myo{\CC}\setminus{G})$ на конденсатор $(E,F)$,
$F=F(f)$. Отметим, что $S$-симметричный компакт $F$ является {\it стационарной
точкой} некоторого функционала энергии (см. п.~\ref{s2s2}).

Обозначим через $\mymu(Q)$ меру, ассоциированную с произвольным полиномом~$Q$:
$$
\mymu(Q)=\sum\limits_{\zeta:Q(\zeta)=0}\delta_\zeta,
$$
где $\delta_\zeta$ --
мера Дирака с носителем в точке $\zeta$. Пусть $\myt\lambda$ -- выметание
равновесной меры $\lambda$ из области Шталя $D=\myo\CC\setminus{F}$
на~$F=\partial{D}$. Тогда в условиях теоремы Шталя для знаменателей
$Q_n(z)$ рациональных функций $R_n$ имеем:
\begin{equation}
\frac1n\mymu(Q_n)\to{\myt\lambda},\qquad n\to\infty,
\label{1.08}
\end{equation}
где сходимость мер понимается в слабой топологии. Равновесная мера $\lambda$
(с носителем на $E$) характеризует предельное распределение точек интерполяции
функции $f$ рациональной функцией~$R_n$.

{\bad
В связи с~\eqref{1.08} (см. также~\eqref{5.31} и приложение~\ref{apl})
напомним хорошо известную теорему Йенча--Сегё (см.~\cite{Jen16},~\cite{Sze21})
для наилучших равномерных {\it полиномиальных} аппроксимаций функции $f$,
голоморфной на $E$.
}

\begin{theoremJS}
Пусть функция $f$ голоморфна на отрезке $E$. Тогда существует бесконечная
подпоследовательность $\Lambda=\Lambda(f)\subset\NN$ такая, что для наилучших
равномерных полиномиальных аппроксимаций $P_n$ функции $f$ имеет место
соотношение:
$$
\frac1n\mu(P_n)\to\lambda_\Gamma,\qquad n\in\Lambda,\quad n\to\infty,
$$
где
$$
d\lambda_\Gamma(\zeta)=
\frac1{2\pi}\frac{|d\zeta|}{|\zeta^2-1|^{1/2}},\quad \zeta\in\Gamma,
$$
-- равновесная мера для эллипса $\Gamma=\Gamma_{\rho_0(f)}$, $\rho_0(f)$ --
индекс голоморфности функции~$f$.
\end{theoremJS}
Таким образом, нули полиномов наилучшего равномерного приближения функции
$f$ в пределе моделируют максимальный эллипс голоморфности~$f$, фактически
``отрезая'' внутренность этого эллипса от его внешней части.
Аналогичный результат справедлив и для наилучших равномерных рациональных
аппроксимаций функции~$f$ с {\it фиксированной} степенью знаменателя.

Отметим, что хорошо известный алгоритм Ремеза~\cite{Rem34} (см.
также~\cite{Akh65}, \cite{Rem69},~\cite{Leb08}) позволяет практически
находить наилучшие равномерные полиномиальные и рациональные аппроксимации {\it
заданной} функции. Этот алгоритм реализован, например, в системе Maple.
Однако процесс практического приближенного построения такой рациональной
аппроксимации предполагает, что исходная функция~$f$ задана в виде явного
аналитического (формульного) выражения.

\subsection{}\label{ss2}
\begin{figure}[h!]
\centerline{
\includegraphics[width=10cm,height=10cm]{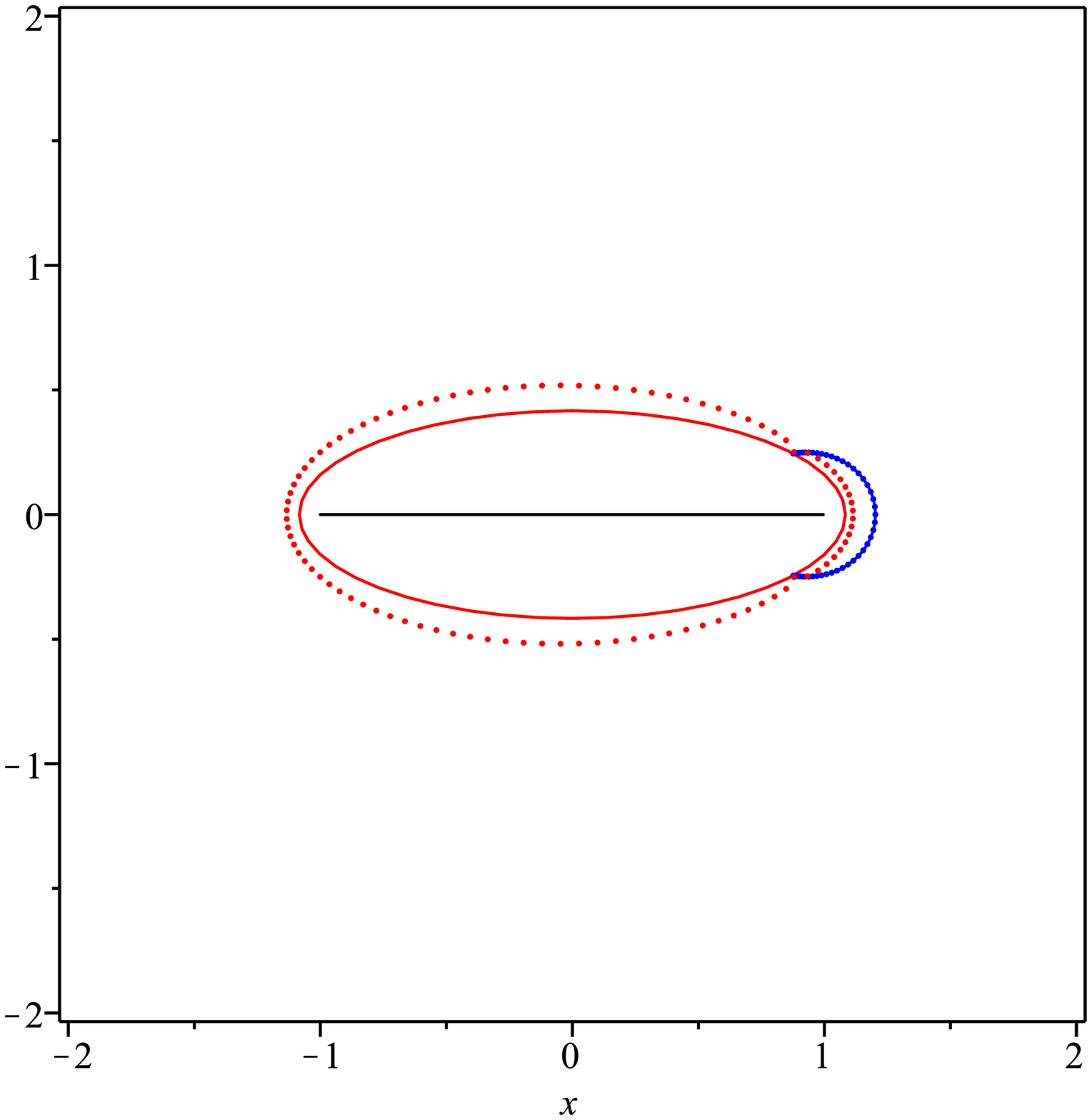}}
\vskip-6mm
\caption{Максимальный эллипс голоморфности (красная линия) функции
$f(z)=\sqrt{(z-a)(z-\myo{a})}$, $\Im{a}>0$, и расположение нулей (красные
точки) частных сумм Фурье--Чебышёва $S_{100}$ и полюсов и нулей (синие точки)
нелинейных аппроксимаций Паде--Чебышёва $F_{50}$.}
\label{Fig1}
\end{figure}

Будем теперь считать, что вещественнозначная функция $f$ {\it задана} на~$E$
(сходящимся) рядом Фурье--Чебышёва (см.~\cite{Boy01}):
\begin{equation}
f(x)=\sum_{k=0}^\infty c_kT_k(x),\qquad c_k=c_k(f)\in\RR.
\label{2.01}
\end{equation}
Отметим, что для частичных сумм $S_n$ ряда~\eqref{2.01}
также справедлив аналог теоремы Йенча--Сегё:
$$
\frac1n\mu(S_n)\to\lambda_\Gamma,
\qquad n\in\Lambda'=\Lambda'(f),\quad n\to\infty.
$$
Таким образом, нули частичных сумм ряда~\eqref{2.01} в пределе также
моделируют максимальный эллипс голоморфности~$f$ (см. рис.~\ref{Fig1}).

Рассмотрим следующую задачу. В классе $\SR_n$ найти рациональную функцию~$F_n$,
голоморфную на отрезке $E$ и удовлетворяющую следующему условию:
\begin{equation}
c_k(F_n)=c_k(f),\qquad k=0,\dots,2n.
\label{2.02}
\end{equation}
{\bad
Иными словами, разложение рациональной функции $F_n=P/Q$ в ряд Фурье--Чебышёва
должно иметь вид
}
\begin{equation}
F_n(x)=c_0+c_1T_1(x)+\dots+c_{2n}T_{2n}+\dotsb.
\label{2.03}
\end{equation}
Подлежат определению из системы~\eqref{2.02} коэффициенты многочленов~$P$
и~$Q$. Рациональная функция $F_n$ (если она существует) называется {\it
нелинейной} диагональной аппроксимацией Паде--Чебышёва функции
$f$ (аппроксимацией Паде ряда~\eqref{2.01}).

{\bad
Система~\eqref{2.02} нелинейна относительно коэффициентов полиномов $P$
и $Q$ и не всегда имеет решение. Тем самым, нелинейная аппроксимация
Паде--Чебышёва может не существовать.
Так как полиномы Чебышёва
являются полиномами Фабера для отрезка $\EE$, то существование
нелинейной аппроксимации Паде--Чебышёва тесно связано с существованием
аппроксимации Паде степенного ряда
$\sum\limits_{k=0}^\infty c_kw^k$, $c_k=c_k(f)$, обладающей определенными
свойствами (см. прежде всего~\cite{Sue80}, а
также~\cite{Ged81},~\cite{Sue09b} и~\cite{Kniz09}).
Отметим, что приведенное определение аппроксимации Паде--Чебышёва (далее --
АЧП) основано на нелинейной (относительно коэффициентов искомой рациональной
функции) схеме Бейкера определения классических аппроксимаций Паде степенного
ряда. Второй способ определения АПЧ, основанный на линейной схеме Фробениуса,
приводит к существенно другим результатам (см.~\cite[часть~2, \S~1.6]{BaGr86}, а
также~\cite{GRS91}, \cite{GRS92},~\cite{GRS10}); сходимости
{\it линейных} аппроксимаций Паде--Чебышёва будет посвящена
работа~\cite{GRS11}. В настоящей работе изучается сходимость {\it нелинейных}
диагональных (типа $(n,n)$) АПЧ.
}

\subsection{}\label{ss3}
В связи с приведенной выше в п.~\ref{ss1} теоремой Шталя и
соотношением~\eqref{1.08} (ср.~\eqref{5.31} и приложение~\ref{apl})
естественным образом возникает вопрос о сходимости нелинейных
АПЧ для такого же класса аналитических функций, имеющих конечное число особых
точек многозначного характера. Настоящая работа посвящена доказательсту аналога
теоремы Шталя для нелинейных АПЧ. Сделаем в этой связи несколько
замечаний.

\subsubsection{}\label{sss1}
В работах авторов~\cite{GRS91} и~\cite{GRS92} для общих ортогональных
разложений были рассмотрены оба способа (Бейкера и Фробениуса) построения
диагональных аппроксимаций Паде и доказаны теоремы о сходимости таких
рациональных аппроксимаций для произвольной марковской функции~\eqref{1.04}
(см. ниже теорему~A). Настоящая работа является естественным развитием
работ~\cite{GRS91} и~\cite{GRS92} для случая нелинейных АПЧ: вместо марковской
функции мы рассматриваем здесь произвольную функцию, вещественную и голоморфную
на~$E$ и имеющую вне $E$ конечное число особых точек многозначного характера.
Основной результат настоящей работы -- теорема~\ref{t2} о скорости сходимости (по
емкости) диагональных АПЧ для функций из указанного класса и о предельном
распределении точек интерполяции функции $f$ рациональной функцией $F_n$. Эта
теорема является аналогом первой теоремы Шталя~\cite{Sta86a} (см. также
приложение~\ref{apl}
к настоящей работе) о сходимости по емкости диагональных аппроксимаций Паде для
функций, заданных (сходящимся) степенным рядом в точке $w=0$ и имеющих в
$\myo{\CC}_w$ конечное число особых точек многозначного характера (ср. теоремой
Шталя для наилучших диагональных аппроксимаций).
При доказательстве сходимости диагональных нелинейных АПЧ мы {\it используем
эту первую теорему Шталя}: благодаря результатам работы~\cite{Sue80} о
нелинейных аппроксимациях Паде--Фабера, нелинейные АПЧ функции $f$ оказываются связанными
с аппроксимациями Паде степенного ряда
$\sum\limits_{k=0}^\infty c_kw^k$, $c_k=c_k(f)$, с помощью оператора
Фабера--Чебышёва (подробнее см. п.~\ref{s4s1}--\ref{s4s2}). Тем самым
{\it непосредственно из первой теоремы Шталя вытекает}, что рациональные
функции~$F_n$ сходятся по емкости к функции $f$ в дополнении к образу компакта Шталя
при отображении, задаваемом функцией Жуковского $z=\Zh(w):=(w+1/w)/2$. Таким образом,
основные доказанные здесь утверждения теоремы~\ref{t2} состоят в следующем:

1) компакт $F$, обладающий $S$-свойством~\eqref{spro1}, существует, единствен
и является стационарным компактом для функционала энергии~\eqref{2.4} при
$\theta=1$;
компакт $F$ совпадает с образом первого компакта Шталя~$S$, заданным в
плоскости переменного $w$, при отображении с помощью функции Жуковского:
$F=\Zh(S)$;

2) справедливо соотношение~\eqref{5.3}, характеризующее скорость сходимости
рациональных аппроксимаций $F_n$ к функции $f$ на компактных подмножествах
области $\myo{\CC}\setminus{F}$ в терминах гринова (относительно стационарного
компакта $F$) потенциала равновесной меры~$\lambda$;
таким образом, для нелинейных АПЧ имеет место аналог формулы~\eqref{1.06},
но с другим $S$-симметричным компактом;

3) предельное распределение точек интерполяции исходной функции $f$
рациональной функцией $F_n$ вполне характеризуется мерой $\lambda$ --
равновесной мерой для стационарного компакта~$F$, а предельное распределение
полюсов $F_n$ -- мерой $\myt{\lambda}$ -- выметанием $\lambda$ из области
$\myo{\CC}\setminus{F}$ на $F$.

\subsubsection{}\label{sss2}
В~\cite{Sue80} (см. также~\cite{Ged81} и~\cite{Sue09b}) был изучен вопрос
о существовании диагональных нелинейных аппроксимаций Паде--Фабера, частным
случаем которых являются аппроксимации Паде--Чебышёва. Оказалось, что при
заданном номере $n$ существование такой аппроксимации тесно связано с
существованием аппроксимации Паде степенного ряда
$\sum\limits_{k=0}^\infty c_kw^k$, $c_k=c_k(f)$, в смысле Бейкера и отсутствием
у этой рациональной функции полюсов в замкнутом единичном круге $|w|\leq1$.
Эти условия выполняются заведомо не при любом $n$. Более того, в~\cite{Sue09b}
приведен пример (основанный на контрпримере Буслаева~\cite{Bus01},~\cite{Bus02}
к гипотезе Бейкера--Гаммеля--Уиллса) гиперэллиптической функции рода $g=2$,
для которой нелинейные АПЧ не существуют ни при одном $n\geq2$. Однако
контрпример Буслаева связан с некоторым ``вырождением'': соответствующая
гиперэллиптическая функция разлагается в периодическую непрерывную дробь.
С другой стороны из работ~\cite{Sta96} и~\cite{Sue00} вытекает, что в
``типичной'' ситуации для достаточно широкого класса аналитических функций,
включающего в себя все ``типичные'' гиперэллиптические функции, всегда существует
некоторая бесконечная подпоследовательность (зависящая от выбранной функции),
для
которой требуемые условия выполняются. Тем самым для соответствующего класса
многозначных функций, полученных с помощью преобразования Жуковского $z=\Zh(w)$,
нелинейные АПЧ существуют по-крайней мере по некоторой подпоследовательности.

\subsubsection{}\label{sss3}
Через $M_1(\EE)$ обозначим множество всех единичных (положительных борелевских)
мер, носители которых принадлежат $\EE$. Пусть $K$ -- произвольный компакт со
связным дополнением в $\myo{\CC}$ такой, что $K\cap\EE=\varnothing$ и область
$D_K=\myo{\CC}\setminus{K}$ регулярна относительно решения задачи Дирихле,
$g_{K}(z,t)$ -- соответствующая области $D_K$ функция Грина с особенностью в
точке~$z=t\in D_K$. Для меры $\mu\in M_1(\EE)$ определены логарифмический
$\displaystyle V^{\mu}(z)=\int_{-1}^1\log|z-x|^{-1}\,d\mu(x)$ и гринов (по
отношению к компакту $K$) потенциалы:
$\displaystyle G^{\mu}_{K}(z)=\int_{-1}^1g_{K}(z,x)\,d\mu(x)$
(полагаем $g_K(z,x)\equiv0$ при $z\in K$, $x\in\EE$). Пусть $\theta\geq0$ --
произвольное фиксированное число. Для фиксированного компакта $K$ существует
{\it единственная} мера $\lambda(\theta)=\lambda_K(\theta)\in\mM_1(\EE)$,
минимизирующая функционал энергии
$J(K,\mu;\theta)
=\displaystyle\int\bigl(\theta V^\mu(x)+G^\mu_{K}(x)\bigr)\,d\mu(x)$
в классе всех мер $\mu\in\mM_1(\EE)$. Мера $\lambda(\theta)$ и только эта мера (в
классе $\mM_1(\EE)$) является {\it равновесной мерой} для смешанного
(гриново-логарифмического) потенциала
$\theta V^\mu(z)+G^\mu_K(z)$.
Другими словами, мера $\lambda(\theta)$ -- единственная мера из класса
$\mM_1(\EE)$, для которой имеет место соотношение равновесия
$\theta V^{\lambda(\theta)}(x)+G^{\lambda(\theta)}_{K}(x)\equiv
w(\theta)=\const$, $x\in \EE$,
$w(\theta)=w_K(\theta)$ -- соответствующая {\it постоянная равновесия}; при
этом $J(K,\mu;\theta)=w(\theta)$.

Теперь для фиксированной функции $f\in\FF(E)$
(см. определения в п.~\ref{s2s3}) и произвольного параметра
$\theta\in[0,+\infty)$ в классе $\KK_f$ допустимых компактов $K$ для $f$
рассмотрим следующую теоретико-потенциальную задачу:
\begin{equation}
\label{thpo1}
\sup_{K\in\KK_f}\inf_{\mu\in M_1(E)}J(K,\mu;\theta)
=\sup_{K\in\KK_f}J(K,\lambda_K;\theta).
\end{equation}
Если существует компакт $F=F(\theta)\in\KK_f$, для которого достигается
супремум в правой части~\eqref{thpo1}, то этот $F$ называется {\it стационарным
компактом} для задачи~\eqref{thpo1}. Значениям параметра $\theta=0,1,3$
соответствуют различные теоретико-потенциальные задачи~\eqref{thpo1} и, вообще
говоря, различные стационарные компакты $F(0),F(1),F(3)$. Эти компакты
совпадают только в исключительных случаях, например, для марковской функции
$f=\myh{\sigma}$, носитель меры которой -- отрезок вещественной прямой
(см.~\cite{GRS91}, \cite{GRS92},~\cite{GRS10}, а также теорему~A). При этом
компакт $F(0)$ соответствует наилучшим рациональным аппроксимациям функции~$f$,
$F(1)$ -- нелинейным АПЧ, $F(3)$ -- линейным АПЧ, и все три стационарных
(в заданном классе $\KK_f$) компакта являются $S$-симметричными (иначе --
обладают $S$-свойством):
\begin{equation}
\frac{\partial G^{\lambda}_{F}}{\partial n_{+}}(\zeta)
=\frac{\partial G^{\lambda}_{F}}{\partial n_{-}}(\zeta),
\qquad\zeta\in F_0,
\label{spro1}
\end{equation}
где $F=F(\theta)$, $\lambda=\lambda_{F}(\theta)\in M_1(E)$ -- соответствующая равновесная
мера, $F_0$ -- объединение всех открытых дуг, принадлежащих компакту $F$,
$\partial/\partial n_{\pm}$ -- нормальные производные, взятые с противоположных
сторон $F_0$. Это свойство, тем самым, носит вполне универсальный характер.

\subsubsection{}\label{sss4}
Несмотря на проблему, связанную с существованием рациональной функции $F_n$,
в системе Maple реализованы именно нелинейные АПЧ. Такие аппроксимации
активно используются в приложениях (см.,
например,~\cite{Kniz84}, \cite{Erm10}, \cite{Boy09},~\cite{Litv03},~\cite{TrGu83},
~\cite{TrGu85}), они имеют
преимущество перед линейными АПЧ, состоящее в том, что по заданным
коэффициентам Фурье--Чебышёва $c_0,c_1,\dots,c_{2n}$ можно
построить нелинейную АПЧ порядка $(n,n)$, а линейную -- только порядка $(m,m)$,
где $m=[2n/3]$ (подробнее см.~\cite{GRS91}, \cite{GRS92}).
Наилучшие равномерные рациональные аппроксимации {\it заданной}
функции могут быть найдены с помощью хорошо известного алгоритма
Ремеза~\cite{Rem69}. Этот алгоритм также реализован в системе Maple.
Однако процесс практического приближенного построения такой рациональной
аппроксимации предполагает, что исходная функция~$f$ задана в виде явного
аналитического (формульного) выражения. Конечного набора
коэффициентов Фурье--Чебышёва для этого недостаточно
(см. рис.~\ref{Fig2}--\ref{Fig4}).

{\bad
В заключение отметим, в~\cite{GRS11} будет изучена сходимость {\it линейных}
АПЧ. В частности, будет показано, что для функция $f\in\FF(\EE)$ при условии
$\rho_0(f)>\sqrt{2}$ существует единственный $S$-симметричный компакт
$F=F(3)\in\KK_f$.
Ограничение $\rho_0(f)>\sqrt2$ связано со спецификой векторной
теоретико-потенциальной задачи равновесия для
параметра $\theta=3$. А именно,
параметрам $\theta=0$, $\theta=1$ и $\theta=3$
соответствуют существенно разные векторные (размера $2\times2$)
теоретико-потенциальные задачи равновесия. При $\theta=0$ матрица
взаимодействия имеет вид
$A_0=\begin{pmatrix}1&-1\\-1&1\end{pmatrix}$, при $\theta=1$ матрица
$A_1=\begin{pmatrix}2&-1\\-1&1\end{pmatrix}$, при $\theta=3$ матрица
$A_3=\begin{pmatrix}4&-1\\-1&1\end{pmatrix}$.
}

\section{Обозначения и формулировка основных результатов}\label{s2}

\subsection{}\label{s2s1}
Перейдем к точным определениям и обозначениям.

Через $M_1(\EE)$ обозначим множество всех единичных (положительных
борелевских) мер, носители которых принадлежат $\EE$.
Пусть $K$ -- произвольный компакт со связным дополнением в $\myo{\CC}$ такой,
что $K\cap\EE=\varnothing$ и область
$D_K=\myo{\CC}\setminus{K}$ регулярна относительно решения задачи Дирихле,
$g_{K}(z,t)$ -- соответствующая области $D_K$ функция Грина с особенностью в
точке~$z=t\in D_K$.
Для меры $\mu\in M_1(\EE)$ определены логарифмический и гринов (по отношению к
компакту $K$) потенциалы:
$$
V^{\mu}(z)=\int_{-1}^1\log\frac1{|z-x|}\,d\mu(x),\qquad
G^{\mu}_{K}(z)=\int_{-1}^1g_{K}(z,x)\,d\mu(x),\qquad z\notin\EE
$$
(полагаем $g_K(z,x)\equiv0$ при $z\in K$, $x\in\EE$).
Пусть $\theta\geq0$ -- произвольное фиксированное число. Для фиксированного
компакта $K$ существует {\it единственная} мера
$\lambda(\theta)=\lambda_K(\theta)\in\mM_1(\EE)$, минимизирующая
функционал энергии
\begin{equation}
J(K,\mu;\theta)=\iint\(\theta\log\frac1{|x-t|}+g_{K}(x,t)\)\,d\mu(x)\,d\mu(t)
=\int\bigl(\theta V^\mu(x)+G^\mu_{K}(x)\bigr)\,d\mu(x)
\label{2.4}
\end{equation}
в классе всех мер $\mu\in\mM_1(\EE)$. Мера $\lambda(\theta)$ и только эта мера (в
классе $\mM_1(\EE)$) является {\it равновесной мерой} для смешанного
(гриново-логарифмического) потенциала
$\theta V^\mu(z)+G^\mu_K(z)$.
Другими словами, мера $\lambda(\theta)$ -- единственная мера из класса
$\mM_1(\EE)$, для которой имеет место соотношение равновесия
\begin{equation}
\theta V^{\lambda(\theta)}(x)+G^{\lambda(\theta)}_{K}(x)\equiv
w(\theta)=\const,\qquad x\in \EE,
\label{2.5}
\end{equation}
$w(\theta)=w_K(\theta)$ -- соответствующая {\it постоянная равновесия}; при этом
$J(K,\mu;\theta)=w(\theta)$.

\subsection{}\label{s2s2}
В~\cite{GRS91} и~\cite{GRS92} была изучена сходимость нелинейных аппроксимаций
$F_n$ и аппроксимаций Фробениуса $\Phi_n$ для общих
ортогональных разложений {\it марковских} функций
\begin{equation}
\myh{\sigma}(z)=\int_{F}\frac{d\sigma(x)}{z-x},\qquad z\in\myo{\CC}\setminus{F},
\label{3.1}
\end{equation}
где $F=[c,d]\subset\RR\setminus\EE$, $\sigma$ -- положительная
борелевская мера на $F$, $\sigma'=d\sigma/dx>0$ почти всюду (п.в.) на $F$.
Скорость сходимости последовательностей $F_n$ и $\Phi_n$ к функции
$f=\myh\sigma$ в области $D=\myo{\CC}\setminus{[c,d]}$ полностью
характеризуется в терминах равновесной меры
$\lambda(\theta)\in M_1(E)$ соответственно для $\theta=1$ и $\theta=3$ следующим
образом.

\begin{theoremA}[(см.~\cite{GRS91},~\cite{GRS92})]
Если $\sigma'>0$ п.в. на $F=[c,d]\subset\RR\setminus\EE$,
то локально равномерно в области $D\setminus\EE$
\begin{equation}
\lim_{n\to\infty}\bigl|(\myh\sigma-f_n)(z)\bigr|^{1/n}
=e^{-2G_{F}^{\lambda(\theta)}(z)}<1,
\label{3.2}
\end{equation}
где $\theta=1$ для $f_n=F_n$ и $\theta=3$ для $f_n=\Phi_n$.
\end{theoremA}

Для наилучших в равномерной метрике на отрезке $\EE$ рациональных аппроксимаций
$R_n$ функции $\myh{\sigma}$ соотношение~\eqref{3.2} справедливо с $\theta=0$
(см.~\cite{Gon78}).

Обозначим через $\mymu(Q)$ меру, ассоциированную с произвольным полиномом~$Q$:
$\mymu(Q)=\sum\limits_{\zeta:Q(\zeta)=0}\delta_\zeta$, где $\delta_\zeta$ --
мера Дирака с носителем в точке $\zeta$. Пусть $\myt\mu$ -- выметание меры
$\mu$ из области $\myo\CC\setminus{F}$ на~$F$. Тогда в условиях
теоремы~A для знаменателей $Q_n(z;\theta)$, $\theta=1,3,0$, соответствующих
рациональных функций $F_n,\Phi_n,R_n$ имеем:
\begin{equation}
\label{count1}
\frac1n\mymu(Q_n(\cdot;\theta))\to{\myt\lambda}(\theta),\qquad n\to\infty,
\end{equation}
где сходимость мер понимается в слабой топологии.

Отметим, что развитые в~\cite{GRS91},~\cite{GRS92} методы позволяют легко
доказать аналог теоремы~A (с заменой равномерной
сходимости~\eqref{3.2} на сходимость по емкости) и для случая, когда $F$
состоит из нескольких отрезков, а $f=\myh\sigma+r$, где $\myh\sigma$~--
марковская функция~\eqref{3.1}, $r$ -- вещественная рациональная функция,
голоморфная на~$\EE$.

\subsection{}\label{s2s3}
{\bad
Введем определения и обозначения, связанные с классом многозначных
аналитических функций, рассматриваемых в настоящей работе. Напомним, что
мы рассматриваем только функции, вещественные и голоморфные на единичном
отрезке $E=[-1,1]$, которые имеют в дополнении к $E$ конечное число особых
точек многозначного характера.
}

Компакт $K$ со связным дополнением $D_K$ в $\myo{\CC}$ будем называть {\it
допустимым} для заданной многозначной аналитической функции~$f$, если
$K\cap\EE=\varnothing$ и $f$ продолжается с отрезка $\EE$ в область $D_K$ как
однозначная мероморфная функция. Множество всех
допустимых компактов для функции~$f$ обозначим через $\KK_f$.

Через $\FF(\EE)$ обозначим класс функций~$f$, вещественных и голоморфных на
$\EE$ и удовлетворяющих следующим двум условиям:

(1) существует конечное множество различных точек
$\Sigma_f=\{b_1,\dots,b_{m}\}
\subset\myo{\CC}\setminus\EE$,
$m=m(f)\geq2$, такое, что:
$\Sigma_f$ симметрично относительно вещественной прямой;
функция $f$ продолжается (с отрезка $\EE$) как многозначная аналитическая
функция в область $\myo{\CC}\setminus\Sigma_f$; каждая точка $b_j\in\Sigma_f$
является точкой ветвления функции~$f$;

(2) существует по-крайней мере один допустимый компакт $K\in\KK_f$ такой,
что $K$ симметричен относительно вещественной прямой, состоит из конечного
числа кусочно аналитических дуг и на каждой открытой дуге, принадлежащей $K$,
скачок функции~$f$ отличен от тождественного нуля.

\begin{figure}[h!]
\centerline{
\includegraphics[width=7.5cm,height=7.5cm]{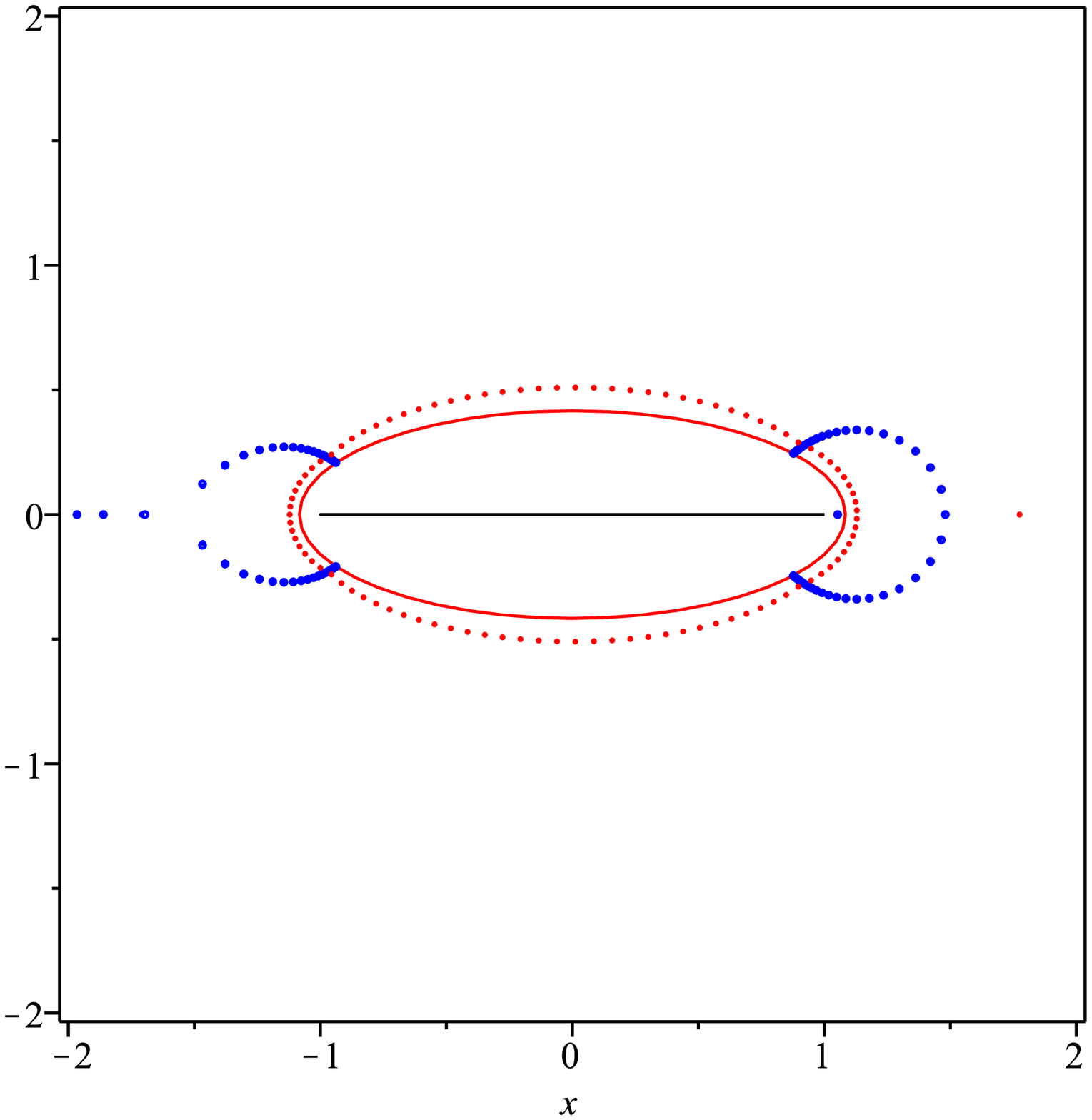}}
\vskip-6mm
\caption{Максимальный эллипс голоморфности (красная линия) функции
$f(z)=\sqrt{(z-a)(z-\myo{a})}+\sqrt[3]{(z-b)(z-\myo{b})(z-c)}$,
$\Im{a},\Im{b}>0$, $c<0$, и расположение нулей (красные
точки) частных сумм Фурье--Чебышёва $S_{100}$ и полюсов и нулей (синие точки)
нелинейных аппроксимаций Паде--Чебышёва $F_{50}$.}
\label{Fig2}
\end{figure}

\begin{figure}[h!]
\centerline{
\includegraphics[width=7.5cm,height=7.5cm]{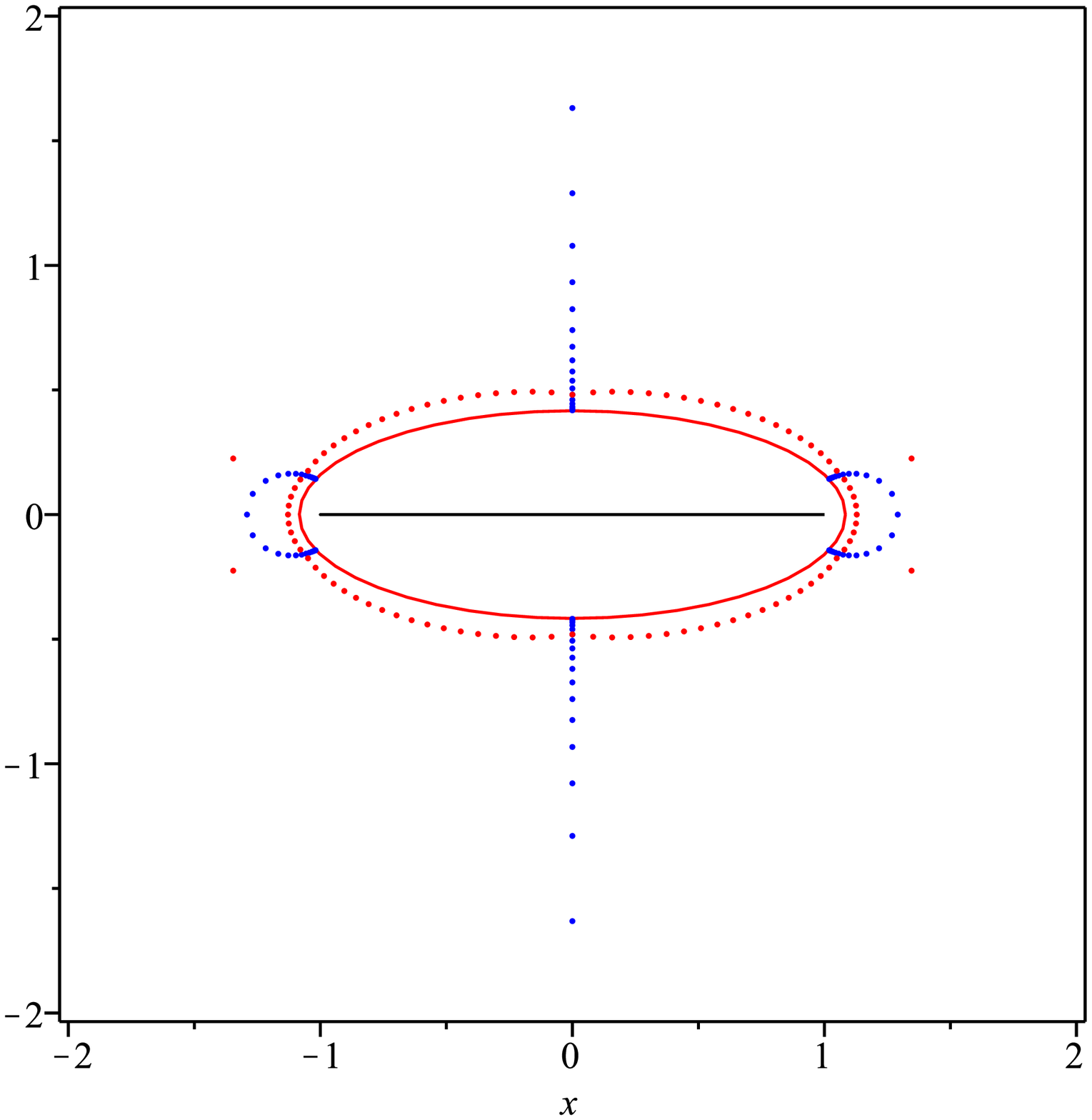}}
\vskip-6mm
\caption{Максимальный эллипс голоморфности (красная линия) функции
$f(z)=\sqrt{(z-a)(z-\myo{a})}+\sqrt{(z-b)(z-\myo{b})}+
\sqrt{(z-ic)(z+ic)}$,
$\Im{a},\Im{b}>0$, $c>0$, и расположение нулей (красные
точки) частных сумм Фурье--Чебышёва $S_{100}$ и полюсов и нулей (синие точки)
нелинейных аппроксимаций Паде--Чебышёва $F_{50}$.}
\label{Fig3}
\end{figure}

\begin{figure}[h!]
\centerline{
\includegraphics[width=7.5cm,height=7.5cm]{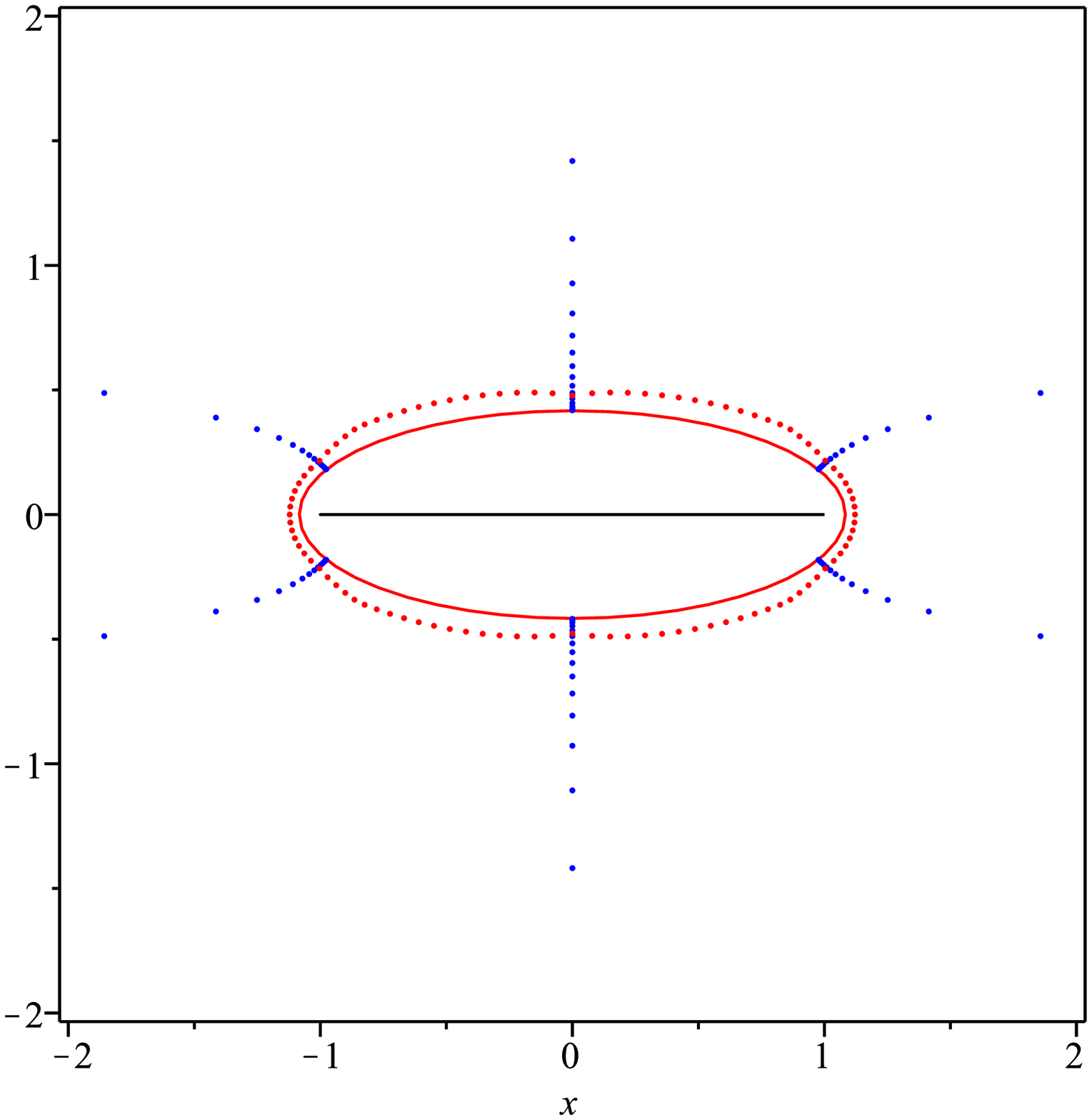}}
\vskip-6mm
\caption{Максимальный эллипс голоморфности (красная линия) функции
$f(z)=\sqrt{(z-a)(z-\myo{a})}+\sqrt{(z-b)(z-\myo{b})}+
\sqrt{(z-ic)(z+ic)}$,
$\Im{a},\Im{b}>0$, $c>0$, и расположение нулей (красные
точки) частных сумм Фурье--Чебышёва $S_{100}$ и полюсов и нулей (синие точки)
нелинейных аппроксимаций Паде--Чебышёва $F_{50}$.}
\label{Fig4}
\end{figure}

Отметим, что классу $\FF(\EE)$ принадлежат, например, следующие функции, не
являющиеся марковскими:
$$
\sqrt{(z-b)(z-\myo{b})},\quad
\root3\of{(z-b)(z-\myo{b})(z-a)},\quad
\log\frac{z-b}{z-\myo{b}},
$$
где $\Im{b}>0$, $a\in\RR\setminus\EE$ и выбрана надлежащая ветвь многозначной
функции (см. рис.~\ref{Fig2}--\ref{Fig4}).

В дальнейшем функция $f\in\FF(E)$ предполагается фиксированной.

\subsection{}\label{s2s4}
Хорошо известно (см.~\cite{Sta85a}--\cite{Sta86b}, \cite{GoRa87}), что при
доказательстве сходимости диагональных аппроксимаций Паде для многозначных
аналитических функций ключевую роль играет существование допустимого (для
заданной функции) компакта, обладающего так называемым свойством {\it
стационарности} для некоторого функционала энергии, соответствующего
рассматриваемой задаче рациональной аппроксимации. Это понятие оказывается
тесно связанным с соответствующей теоретико-потенциальной задачей равновесия,
а стационарный компакт оказывается $S$-симметричным (иначе говоря, обладает
$S$-свойством).
Приведем определение $S$-свойства допустимого компакта, соответствующего
рассматриваемой задаче равновесия~\eqref{2.5}.

\begin{definition}\label{de1}
Пусть параметр $\theta\geq0$. Будем говорить, что
(не разбивающий плоскость и состоящий из конечного числа кусочно-аналитических
дуг) допустимый компакт $F=F(\theta)\in\KK_f$ обладает $S$-{\it свойством}
(или является $S$-симметричным), если:
\begin{equation}
\frac{\partial G^{\lambda}_{F}}{\partial n_{+}}(\zeta)
=\frac{\partial G^{\lambda}_{F}}{\partial n_{-}}(\zeta),
\qquad\zeta\in F_0,
\label{spro2}
\end{equation}
где $\lambda=\lambda_F(\theta)$ -- соответствующая равновесная мера,
$F_0$ -- объединение всех открытых дуг, принадлежащих компакту $F$,
$\partial/\partial n_{\pm}$ -- нормальные производные, взятые с противоположных
сторон $F_0$.
\end{definition}

Зафиксируем теперь параметр $\theta=1$ и в дальнейшем в обозначениях будем как
правило опускать указание на этот параметр. Пусть $\lambda=\lambda_K(1)\in
M_1(E)$
-- равновесная мера, соответствующая произвольному компакту $K\in\KK_f$,
$w=w_K(1)$ -- соответствующая постоянная равновесия (см.~\eqref{2.5}):
$V^{\lambda}(x)+G^{\lambda}_{K}(x)\equiv w$, $x\in \EE$;
при этом $J(K,\lambda)=\min\limits_{\mu\in\mM_1(\EE)}J(K,\mu)=w$.

Справедлива следующая
\begin{theorem}\label{t1}
Если функция $f\in\FF(\EE)$, то существует единственный компакт
$F=F(1)\in\KK_f$ такой,
что
\begin{equation}
J(F,\lambda_{F})=\max_{K\in\KK_f}J(K,\lambda_K).
\label{5.1}
\end{equation}
Стационарный компакт $F$ состоит из конечного числа кусочно-аналитических дуг, не
разбивает плоскость и обладает $S$-свойством~\eqref{spro2},
где $\lambda_F=\lambda_F(1)$ -- соответствующая равновесная мера.
\end{theorem}

Теорема~\ref{t1} доказывается в два этапа в соответствии со следующей схемой.
Сначала с помощью геометрических соображений, основанных
на замене функционала энергии~\eqref{2.4} меры $\lambda\in\mM_1(\EE)$ на
обобщенный (по отношению к допустимому компакту $K\in\KK_f$)
трансфинитный диаметр~$\EE$ доказывается, что максимум в правой
части~\eqref{5.1} достаточно искать среди тех допустимых компактов, которые
лежат вне максимального канонического эллипса голоморфности~$f$. Такое
семейство компактно в хаусдорфовой метрике, поэтому существует допустимый
компакт~$F$, удовлетворяющий соотношению~\eqref{5.1}. Затем с помощью
вариационного метода аналогично~\cite{PeRa94} устанавливается, что этот
экстремальный компакт~$F$ является замыканием критических траекторий некоторого
квадратичного дифференциала. Отсюда уже вытекает $S$-свойство~\eqref{spro2}.

Справедлива следующая теорема о скорости сходимости рациональных аппроксимаций
$F_n$ к функции $f$ в области $D=D_{F}=\myo{\CC}\setminus{F}$.
\begin{theorem}\label{t2}
Пусть $f\in\FF(\EE)$. Тогда для любого компакта
$K\subset\myo{\CC}\setminus(\EE\cup{F})$
\begin{equation}
\bigl|(f-F_n)(z)\bigr|^{1/n}\overset{\mcap}{\longrightarrow}
e^{-2G_{F}^{\lambda_{F}}(z)}<1,
\qquad z\in K,
\label{5.3}
\end{equation}
где $F=F(1)$. При этом для точек интерполяции функции $f$ рациональной функцией
$F_n$ на отрезке $E$ и знаменателя $Q_n$ имеем
\begin{equation}
\label{5.31}
\frac1{2n}{\mu(\omega_{2n})}\to\lambda_{F},\quad
\frac1n\mymu(Q_n)\to\myt\lambda_{F},\qquad n\to\infty,
\end{equation}
где $\omega_{2n}$ полином степени $2n+1$ с нулями в точках интерполяции,
$\myt{\lambda}_{F}$ -- выметание меры $\lambda_{F}\in\mM_1(\EE)$
из области $D$ на $\partial D=F$.
\end{theorem}
Тем самым в условиях теоремы~\ref{t2} последовательность
$F_n\overset{\mcap}\longrightarrow f$ на компактных подмножествах области
$D=\myo{\CC}\setminus{F}$.

Доказательство теоремы~\ref{t2} основано на первой теореме Шталя (см.
приложение~\ref{apl}) и общей схеме, предложенной
в~\cite{Sta85a}--\cite{Sta86b} и \cite{GoRa87} и основанной на
$S$-свойстве~\eqref{spro2} соответствующего стационарного компакта.
При доказательстве предельного соотношения~\eqref{5.31} используется общая
теоремы~3 работы~\cite{GoRa87} (для рассматриваемого здесь частного случая
$\psi(z)=-V^\lambda(z)$). Отметим, что из~\eqref{5.3} вытекает, что каждый
полюс~$f$ в~$D$ притягивает по-крайней мере столько полюсов $f_n$, какова его
кратность.

Таким образом, для функции~$f$ класса $\FF(\EE)$ компакт $F(1)$ в случае
нелинейных АПЧ играет роль отрезка $[c,d]$, соответствующего марковской функции
$\myh\sigma$, $\supp{\sigma}=[c,d]$ (см.~\eqref{3.1}).

Еще раз отметим, что стационарный компакт $F(1)$ является
образом компакта Шталя для степенного ряда $\sum\limits_{k=0}^\infty c_kw^k$,
$c_k=c_k(f)$, при отображении, задаваемом функцией Жуковского $z=(w+w^{-1})/2$.

\subsection{}\label{s2s7}
Если (для параметра $\theta=1$) стационарный компакт $F=F(1)$,
соответствующий функции $f\in\FF(E)$, состоит из $s$ непересекающихся
аналитических дуг, попарно соединяющих некоторые точки ветвления
$b'_1,\dots,b'_{2s}\in\Sigma_f$ функции~$f$, то он допускает наглядное описание
в терминах, связанных с {\it четырехлистной} римановой поверхностью рода
$g=s-1$.

Построим сначала двулистную риманову поверхность
$\RS=\RS^{(1)}\cup\RS^{(2)}$ следующим образом. Возьмем два экземпляра
римановой сферы $\myo{\CC}$, разрезанных по отрезку $\EE=[-1,1]$, и переклеим
по разрезу. Полученная двулистная риманова поверхность
$\RS=\RS^{(1)}\cup\RS^{(2)}$ эквивалентна римановой сфере. Определим на $\RS$
функцию $u(\zz)$, $\zz\in\RS$, следующим образом:
$u(z^{(1)})=G^{\lambda_{F}}_{F}(z)$, $u(z^{(2)})=w_F-V^{\lambda_{F}}(z)$.
Непосредственно из условия равновесия~\eqref{2.5} вытекает, что $u$ --
гармоническая функция на $\RS\setminus(F^{(1)}\setminus\{\infty^{(2)}\})$.
Кроме того, $u\equiv0$ на компакте $F^{(1)}$ и
$u(z^{(2)})=\log|z|+w_F+o(1)$ при $z^{(2)}\to\infty^{(2)}$.
Возьмем теперь два экземпляра $\RS$, разрезанных по компакту $F^{(1)}$, и
переклеим их между собой по
соответствующим разрезам. Получим четырехлистную риманову поверхность $\RS_1$
рода $g=s-1$. Так как $u\equiv0$ на $F^{(1)}$, то $u$
гармонически продолжается с одного экземпляра $\RS$ на другой с переменой
знака. Продолженная функция гармонична на $\RS_1$ всюду кроме точек
$\zz=\infty^{(2)},\infty^{(3)}$, где она имеет логарифмические особенности:
$\log{|z|}$ при $\zz\to\infty^{(2)}$ и
$-\log{|z|}$ при $\zz\to\infty^{(3)}$.
Следовательно, $u(\zz)=\Re\Omega(\zz)$,
$\Omega(\zz)=\displaystyle\int_{b_1}^{\zz}d\Omega$, где $d\Omega(\zz)$ --
(единственный) абелев
дифференциал на $\RS_1$ с чисто мнимыми периодами и особенностью вида
$1/z$ в точке $\zz=\infty^{(2)}$ и вида
$-1/z$ в точке $\zz=\infty^{(3)}$. Компакт $F$
соответствует нулевой линии уровня
функции $\Re\Omega(\zz)$: $F=\{z\in\myo{\CC}:\Re\Omega(\zz)=0\}\setminus\EE$.

Отметим, что $u(\zz)=g_{F^{(1)}}(\zz,\infty^{(2)})$ --
функция Грина для области $\RS\setminus{F^{(1)}}$ с особенностью в точке
$\zz=\infty^{(2)}$, а $w_F=\gamma^{(2)}$ -- постоянная Робена для этой
функции Грина. Тем самым задача о максимуме
постоянной $w_F$ соответствует задаче о минимуме $e^{-\gamma}$, т.е. минимуме
логарифмической емкости.

Непосредственно из результатов
работ~\cite{Sta85a}--\cite{Sta86b} и~\cite{GoRa87} вытекает, что для {\it
наилучших} в равномерной метрике на отрезке $[-1,1]$ рациональных аппроксимаций
$f_n=R_n$ функции $f\in\FF(\EE)$ также справедливо
соотношение вида~\eqref{5.3}, где $\theta=0$, $F=F(0)$ -- стационарный
компакт, соответствующий задаче равновесия~\eqref{2.5} с
$\theta=0$ и обладающий $S$-свойством~\eqref{spro2}, $\lambda=\lambda_F(0)$ --
соответствующая равновесная мера.
В этом случае функция $u(z^{(1)})=w_F-G^\lambda(z)$ продолжается через разрез,
проведенный по отрезку $E$, на второй лист римановой поверхности
$\RS=\RS^{(1)}\cup\RS^{(2)}$ с переменой знака. Отсюда уже легко вытекает,
что задача о максимуме постоянной $w_F$ эквивалентна задаче о минимуме
емкости конденсатора $(F^{(1)},F^{(2)})$.

Таким образом, все три стационарных (в заданном классе $\KK_f$)
компакта $F(1),F(3),F(0)$ обладают $S$-свойством~\eqref{spro2}. Это свойство, тем
самым, носит вполне универсальный характер.

\section{Доказательство теоремы~\ref{t1}}\label{s3}

\subsection{}\label{s3s1}
Пусть $K\in\KK_f$ -- допустимый компакт для функции $f\in\FF(E)$:
$K\subset\myo{\CC}\setminus{E}$, состоит из конечного числа
кусочно-аналитических дуг, не разбивает плоскость,
$f\in\MM(\myo{\CC}\setminus{K})$ и для любой дуги $\ell\subset{K}$ скачок
$\Delta{f}\not\equiv0$ на $\ell$. Предположим также, что $\infty\notin{K}$,
$\infty\notin{\Sigma_f}$.

Пусть $\lambda\in M_1(E)$ -- равновесная мера для $K$:
\begin{equation}
V^\lambda(x)+G^\lambda_K(x)\equiv w_K,\qquad x\in{E},
\label{51.1}
\end{equation}
величина $w=w_K$ -- ``постоянная Робена'' для $K$, $e^{-w}$ -- ``емкость''
компакта~$K$ (соответствующая энергии взаимодействия $J(\cdot)$;
см.~\eqref{2.4}). Пусть $\myt{\lambda}$ выметание равновесной меры $\lambda$ на
$K$. Рассмотрим функцию
\begin{equation}
v(z):=V^\lambda(z)+G_K^\lambda(z)+G_E^{\myt\lambda}(z)+g_E(z,\infty),
\label{51.2}
\end{equation}
где $G_E^\mu(z)$ -- гринов относительно $E$ потенциал меры $\mu$,
$$
G_E^\mu(z)=\int g_E(z,\zeta)\,d\mu(\zeta),\qquad
\supp{\mu}\subset\myo{\CC}\setminus{E},
$$
$g_E(z,\zeta)$ -- функция Грина для дополнения к отрезку $E$,
$g_E(z,\infty)=\log|z+\sqrt{z^2-1}|$. Нетрудно видеть, что сумма
$G_K^\lambda(z)+G_E^{\myt\lambda}(z)$ -- функция, гармоническая вне $E$.
Из~\eqref{51.2} получаем, что и функция $v(z)$ гармонична в
$\myo{\CC}\setminus{E}$. Из соотношения равновесия~\eqref{51.1} вытекает, что
$v(x)\equiv w_K$ при $x\in{E}$. Следовательно, $v(z)\equiv w_K$ при
$z\in\myo{\CC}$. Отсюда для $z\in{K}$ получаем
\begin{equation}
w_K\equiv v(z)=V^\lambda(z)+G^{\myt{\lambda}}_E(z)+g_E(z,\infty),
\qquad z\in K.
\label{51.3}
\end{equation}
Но при $z\in{K}$ потенциал $V^\lambda(z)=V^{\myt{\lambda}}(z)-c_K$, где
$c_K=\displaystyle\int_E g_K(x,\infty)\,d\lambda(x)$.
Таким образом, из~\eqref{51.3} получаем, что для выметания $\myt\lambda\in
M_1(K)$ равновесной меры
$\lambda\in M_1(E)$ справедливо равенство:
\begin{equation}
V^{\myt{\lambda}}(z)+G^{\myt{\lambda}}_E+g_E(z,\infty)\equiv\const=w_K+c_K=:w_E,
\qquad z\in{K}.
\label{51.4}
\end{equation}
Следовательно, мера $\myt\lambda\in M_1(K)$ является (единственной) равновесной
мерой в классе $M_1(K)$ для потенциала $V^\mu(z)+G^\mu_E(z)$ с внешним полем
$\pfi(z)=g_{E}(z,\infty)$.
Таким образом, отображение $\lambda\to\myt\lambda$ задает взаимно однозначное
соответствие между {\it равновесными} мерами для задач~\eqref{51.1}
и~\eqref{51.4}.
Хорошо известно (см.~\cite{Lan66},~\cite{GoRa81},~\cite{GoRa84},~\cite{GoRa87}), что
равновесная мера и только эта мера в классе мер $\mu\in M_1(K)$ минимизирует
соответствующий функционал энергии
\begin{align}
J_\pfi(K,\mu):&=\iint\biggl\{\log\frac1{|z-\zeta|}+g_E(z,\zeta)\biggr\}
\,d\mu(z)d\mu(\zeta)+2\int\pfi(z)\,d\mu(z)\notag\\
&=\int\biggl\{V^\mu(z)+G_E^\mu(z)+\pfi(z)\biggr\}\,d\mu(z)
+\int\pfi(z)\,d\mu(z)=w_E+c_E.
\label{51.5}
\end{align}
При этом для соответствующих задачам~\eqref{51.1} и~\eqref{51.4} постоянных
равновесия имеем: $w_K=w_E-c_K$. Так как $v(z)\equiv w_K$, то из~\eqref{51.3}
получаем, что $w_K=v(\infty)=\gamma_E+c_K+c_E$, где
$$
c_K=\int_{E}g_K(x,\infty)\,d\lambda(x),\qquad
c_E=\int_{K}g_E(\zeta,\infty)\,d\myt\lambda(\zeta),
$$
$\gamma_E=\log2$ -- постоянная Робена для $E$. Отсюда с учетом~\eqref{51.5} и
равенства
$w_K=\gamma_E+c_K+c_E$ уже легко получаем, что
соответствующие равновесным мерам энергии связаны соотношением:
\begin{equation}
J_\pfi(K,\myt\lambda)=2J(K,\lambda)-\gamma_E.
\label{51.7}
\end{equation}
Непосредственно из~\eqref{51.7} вытекает, что задача
\begin{equation}
\sup_{K\in\KK_f}\inf_{\nu\in M_1(E)}J(K,\nu)
\label{51.8}
\end{equation}
эквивалентна (при условии, что решение ищется среди допустимых
компактов $K\not\ni\infty$) задаче
\begin{equation}
M=\sup_{K\in\KK_f}\inf_{\mu\in M_1(K)}J_\pfi(K,\mu),\quad\text{ где }
\quad\pfi(z)=g_{E}(z,\infty).
\label{51.9}
\end{equation}

Ниже мы докажем {\it существование\/} стационарного компакта $F$ для
задачи~\eqref{51.9}. В силу указанной эквивалентности это будет означать,
что этот компакт $F$ является стационарным компактом и для
задачи~\eqref{51.8} при условии, что $F\not\ni\infty$. Единственность такого
компакта будет вытекать как обычно из сходимости АПЧ к функции $f\in\MM(D_F)$.
Условия $F\not\ni\infty$ и $\Sigma\not\ni\infty$ не
являются ограничительными и носят чисто технический характер: их всегда можно
добиться с помощью подходящего дробно линейного преобразования, сохраняющего
отрезок $E$ (пока речь идет не об аппроксимации ряда~\eqref{1.01} рациональными
функциями, а о теоретико-потенциальных задачах
равновесия~\eqref{51.8},~\eqref{51.9}). Отметим, что непосредственно для
исследования сходимости аппроксимаций нужно
эквивалентное стационарности свойство
$S$-симметрии компакта $F$. Как следует из определения
(см.~\eqref{spro2}), это свойство инвариантно относительно дробно-линейных
преобразований и вообще не зависит от того, содержит компакт $F$ точку
$z=\infty$ или нет.

В заключение этого пункта приведем еще одну эквивалентную переформулировку
экстремальной скалярной задачи~\eqref{51.8} в виде матричной (размера
$2\times2$) теоретико-потенциальной задачи равновесия
(см.~\cite{GoRa87},~\cite{GoRa80}). Рассмотрим следующую
векторную задачу равновесия в классе вектор-мер $\vec{\mu}=(\mu_1,\mu_2)$,
$\mu_1\in M_1(E),\mu_2\in M_1(K)$, $K\in\KK_f$:
\begin{equation}
\begin{aligned}
2V^{\mu_1}(x)-V^{\mu_2}(x)&\equiv w_1,\quad x\in E,\\
-V^{\mu_1}(z)+V^{\mu_2}(z)&\equiv w_2,\quad z\in K,\\
\end{aligned}
\label{51.10}
\end{equation}
с матрицей взаимодействия $A=\begin{pmatrix} 2&-1\\-1&1\end{pmatrix}$.
Векторная задача~\eqref{51.10} эквивалентна скалярной задаче
\begin{equation}
V^\mu(x)+G_K^\mu(x)\equiv w,\qquad x\in E,\quad\mu\in M_1(E).
\label{51.11}
\end{equation}
При этом решение $\vec{\lambda}=(\lambda_1,\lambda_2)$ задачи~\eqref{51.10} связано с решением
задачи~\eqref{51.11} следующим образом: $\lambda_1=\lambda$,
$\lambda_2=\myt\lambda$ -- выметание меры $\lambda\in M_1(E)$ из
$\myo{\CC}\setminus{K}$ на $K$, а для соответствующих постоянных равновесия
$w_1,w_2$ и $w$ имеем: $w_2=c_K(\lambda)=\int_{E}g_K(x,\infty)\,d\lambda(x)$,
$w_1=w-w_2$; тем самым, $w=w_1+w_2$. Таким образом, скалярная экстремальная
задача~\eqref{51.8} может быть эквивалентным образом переформулирована в
виде следующей векторной задачи
\begin{equation}
\sup_{K\in\KK_f}\inf_{\substack{(\mu_1,\mu_2),\\\mu_1\in M_1(E),\mu_2\in M_1(K)}}
J_A(K,\vec{\mu}),
\label{51.12}
\end{equation}
где
$$
J_A(K,\vec{\mu}):=(A\vec{\mu},\vec{\mu})=\sum_{k,j=1}^2 a_{k,j}[\mu_k,\mu_j],
\quad
[\mu_k,\mu_j]=\iint_{K\times K}\log\frac{1}{|z-\zeta|}\,d\mu_k(z)d\mu(\zeta).
$$

\subsection{}\label{s3s2}
Следующие п.~\ref{s3s2}--\ref{s3s4} посвящены доказательству существования
стационарного компакта для экстремальной задачи~\eqref{51.9} и описанию
некоторых его свойств, в том числе $S$-свойства~\eqref{spro2}. Схема изложения следующая. Сначала
мы докажем, что решение экстремальной задачи~\eqref{51.9} (т.е.
стационарный компакт $F$) можно искать среди тех допустимых компактов, которые
лежат вне максимального эллипса голоморфности функции $f$. Затем пользуясь
тем, что такое подсемейство компактно в хаусдорфовой топологии мы докажем
существование допустимого компакта, на котором достигается
величина~\eqref{51.9}. Затем с помощью вариационного метода
аналогично~\cite{PeRa94} доказывается, что экстремальный компакт является
замыканием критических траекторий некоторого квадратичного дифференциала.
Отсюда уже вытекает $S$-свойство~\eqref{spro2}.
Отметим, что существует и другой подход к задаче о существовании и свойствах
стационарного компакта. Этот подход основан на явном описании такого компакта
в терминах некоторой алгебраической функции (см.,
например,~\cite{ApLyTu11},~\cite{ApKuVa08}).

В классе мер $\mu\in M_1(K)$ рассмотрим следующую задачу равновесия:
\begin{equation}
V^{\mu}(z)+G^{\mu}_E(z)+g_E(z,\infty)\equiv\const=\myt{w},
\qquad z\in{K},
\label{61.1}
\end{equation}
Существует единственная {\it равновесная мера\/} $\myt\lambda\in M_1(K)$, для
которой справедливо~\eqref{61.1} с некоторой постоянной $\myt w=\myt w_K$
(см.~\cite{Lan66}, а также \cite{GoRa81},~\cite{GoRa84},~\cite{GoRa87}).
Хорошо известно (см.~\cite{GoRa81},~\cite{GoRa84},~\cite{GoRa87}), что
равновесная мера и только эта мера минимизирует соответствующий функционал
энергии
\begin{align}
J_\pfi(K,\mu):&=\iint\biggl\{\log\frac1{|z-t|}+g_E(z,t)\biggr\}
\,d\mu(z)d\mu(t)+2\int\pfi(z)\,d\mu(z)\notag\\
&=\int\biggl\{V^\mu(z)+G_E^\mu(z)+\pfi(z)\biggr\}\,d\mu(z)
+\int\pfi(z)\,d\mu(z)\notag\\
&=\iint\biggl\{\log\frac1{|z-t|}+g_E(z,t)+\pfi(z)+\pfi(t)\biggr\}
\,d\mu(z)d\mu(t),
\label{61.2}
\end{align}
где $\pfi(z)=g_E(z,\infty)=\log|z+\sqrt{z^2-1}|$.

Рассмотрим следующую задачу:
\begin{equation}
M=\sup_{K\in\KK_f}\min_{\mu\in M_1(K)}J_\pfi(K,\mu),\quad\text{ где }
\quad\pfi(z)=g_{E}(z,\infty).
\label{61.3}
\end{equation}
Докажем, что существует локальное решение задачи~\eqref{61.3}, принадлежащее
{\it внешности} максимального эллипса голоморфности функции $f$. Для этого
определим внешнее проектирование $\widehat{K}$ произвольного компакта
$K\in\KK_f$ вдоль гипербол, соответствующих отрезку $E=[-1,1]$, следующим
образом. Если точка $z\in\Omega=\myo\CC\setminus{D_\rho}$, то $\widehat z=z$.
Пусть теперь $z\in{D_\rho}\setminus{E}$. Тогда через эту точку проходит
единственная гипербола с фокусом в точке $+1$ или $-1$; так как $z$ не лежит на
отрезке $E$, то она принадлежит верхней или нижней ветви этой гиперболы. Положим
$\widehat z=z_0\in\Gamma_\rho$ -- единственная точка пересечения этой ветви с
максимальным эллипсом голоморфности функции~$f$. Ясно, что при таком
проектировании $\widehat K\in\KK_f$ для любого $K\in\KK_f$.

\begin{lemma}\label{l1}
Для любого компакта $K\in\KK_f$ имеем:
$J_\pfi(\widehat K,\myt\lambda_{\widehat K})\geq J_\pfi(K,\myt\lambda_{K})$.
\end{lemma}

\begin{proof}
Для $z\in\myo\CC\setminus{E}$ положим $z=(\zeta+1/\zeta)/2$, где $|\zeta|<1$.
Тогда введенная выше операция внешнего проектирования вдоль
гипербол в плоскости $\CC_z$ эквивалентна радиальному проектированию внутрь
единичного круга $|\zeta|<1$ в плоскости $\CC_\zeta$, а для функций Грина
отрезка $E$ имеем:
\begin{equation}
g_{E}(z,\infty)=\log\frac1{|\zeta|},\qquad
g_E(z,t)=\log\biggl|\frac{1-\myo{\xi}\zeta}{\zeta-\xi}\biggr|,
\label{61.4}
\end{equation}
где $t=(\xi+1/\xi)/2$, $|\xi|<1$.
Рассмотрим теперь аналог функционала~\eqref{61.2} в классе дискретных мер
$\mu\in M_1(K)$, $K\in\KK_f$ (т.е. аналог трансфинитного диаметра для
допустимых компактов; подробнее см.~\cite[глава~II, \S\,3]{Lan66}). Так как по условию компакт $K$ симметричен
относительно вещественной оси, то можно рассматривать только меры с
симметричным относительно вещественной оси носителем. Соответствующий класс
дискретных мер обозначим $M_1^{*}(K)$. Тогда
$$
\inf_{\mu\in M_1(K)}I_\pfi(K,\mu)=
\inf_{\mu\in M_1^{*}(K)}J_\pfi(K,\mu).
$$
С учетом~\eqref{61.4} подынтегральное выражение в~\eqref{61.2} переписывается
следующим образом:
\begin{align}
\log\frac1{|z-t|}
&+g_E(z,t)+\pfi(z)+\pfi(t)
=\log\frac{e^{g(z,\infty)}e^{g(\zeta,\infty)}}{|z-\zeta|}+g_E(z,\zeta)\notag\\
&=\log(|\zeta-\xi|\cdot|1-\zeta\xi|)
+\log\frac{|\zeta-\xi|}{|1-\myo{\xi}\zeta|}.
\label{61.5}
\end{align}
Для дискретной меры с симметричным носителем функционал энергии~\eqref{61.2}
не изменится, если мы в~\eqref{61.5} уберем знак комплексного сопряжения и
полученное выражение $2\log|\zeta-\xi|$ подставим в~\eqref{61.2}. Теперь
уже ясно, что величина $\inf\limits_{\mu\in M_1^{*}(K)}I_\pfi(K,\mu)$
соответствует емкости прообраза компакта $K$ при отображении $z=\Zh(\zeta)$,
$|\zeta|<1$. Хорошо известно, что при радиальном проектировании емкость
множества не увеличивается. Лемма~\ref{l1} доказана.
\end{proof}

Итак, мы показали, что решение экстремальной задачи~\eqref{51.9} (т.е.
стационарный компакт $F$) можно искать среди тех допустимых компактов, которые
лежат вне максимального эллипса голоморфности функции $f$; более точно, в таком
семействе компактов (подсемействе $\KK_f$) мы будем искать {\it локальное}
решение экстремальной задачи~\eqref{51.9} при априорном условии, что
соответствующий компакт $F\not\ni\infty$. Точнее, зафиксируем окрестность
$U(\infty)$ бесконечно удаленной точки $z=\infty$ и в дальнейшем считаем, что
супремум можно брать по подсемейству тех компактов из $\KK_f$, которые не
пересекаются с $U(\infty)$.

Величина $M<\infty$. Следовательно, существует последовательность
допустимых компактов $K_n\in\KK_f$ такая, что все
$K_n\subset\CC\setminus{D_{\rho_0(f)}}$ и
$J_\pfi(K_n,\myt\lambda_{K_n})\to{M}$.
Семейство $\{K_n\}$ компактно в хаусдорфовой топологии. Значит, существует
подпоследовательность, которую мы также будем обозначать $\{K_n\}$,
и компакт $K^{*}\subset\CC\setminus{D_{\rho_0(f)}}$ такие, что
$d_{\rH}(K_n,K^{*})\to0$ при $n\to\infty$,
$K^{*}\subset\CC\setminus{D_{\rho_0(f)}}$ (здесь и далее через
$d_{\rH}(\cdot,\cdot)$ обозначается расстояние между двумя компактами в
хаусдорфовой метрике).

Так как функция $f$ имеет конечное число точек ветвления, то для всех
достаточно больших $n\geq n_0$ все компакты $K_n$ содержат в точности одни
и те же точки ветвления функции $f$. Значит, это справедливо и для компакта
$K^{*}$.
Очевидно, что $K^{*}$ симметричен относительно вещественной прямой и функция
$f$ продолжается как (однозначная) мероморфная функция в ту связную компоненту
$D^{*}$ дополнения к $K^{*}$, которая содержит отрезок $E$. Компакт $K^{*}$
состоит из конечного числа континуумов, поэтому область $D^{*}$ регулярна
относительно решения задачи Дирихле. Следовательно, существует функция Грина
$g_{K^{*}}(z,\zeta)$ для области $D^{*}$.

Приведем некоторые свойства компакта $K^{*}$.

Компакт $K^{*}$ не имеет внутренних точек и не разбивает плоскость.
Доказательство этого утверждения проводится по схеме, предложенной
в~\cite{BerRul02} (см. также~\cite{MarRak09a},~\cite{MarRak09b}).
Компакт $K^{*}$ состоит из конечного числа
континуумов, каждая связная компонента $K^{*}$ содержит по-крайней мере две
точки ветвления функции~$f$. Так как $d_{\rH}(K_n,K^{*})\to0$ при $n\to\infty$,
то для
$$
\eps=\frac14\min_{k\neq j,b_j,b_k\in\Sigma_f}|b_j-b_k|
$$
имеем: $d_{\rH}(K_n,K^{*})<\eps$, $d_{\rH}(K_n,K_m)<\eps$ при $n,m\geq n_0$;
следовательно, при $n\geq n_0$ число связных компонент $K_n$ и $K^{*}$
одинаково, сходимость $K_n\to K^{*}$ можно трактовать по-компонентно и
соответствующие компоненты $K_n$ и $K^{*}$ содержат одни и те же точки
множества $\Sigma_f$; функция $f\in\MM(D_{E}(K^{*}))$, граница области
$D_{E}(K^{*})$ состоит из конечного числа континуумов и $D_{E}(K^{*})$
регулярна относительно решения задачи Дирихле; существует функция Грина
$g_{K^{*}}(z,\zeta)$ для области $D_{E}(K^{*})$.

По условию выбора последовательности $\{K_n\}$ имеем:
$J_\pfi(K_n,\myt\lambda)\to M$. В силу определения величины $M$ (см.~\eqref{51.9}) справедливо неравенство
\begin{equation}
J_\pfi(K^{*},\myt\lambda_{K^{*}})\leq{M}.
\label{51.12}
\end{equation}
Так как для равновесных мер $\myt\lambda_K$ функционал $J(K,\lambda_K)$
(см.~\eqref{2.4}) полунепрерывен сверху по $K$ и
$J_\pfi(K,\myt\lambda)=2J(K,\lambda)-\gamma_E$, то из соотношения
$J_\pfi(K_n,\myt\lambda)\to M$ вытекает, что
$$
M=\varlimsup_{n\to\infty}J_\pfi(K_n,\myt\lambda)\leq J_\pfi(K^{*},\myt\lambda_{K^{*}}).
$$
Отсюда и неравенства~\eqref{51.12}
получаем, что $J_\pfi(K^{*},\myt\lambda_{K^{*}})=M$.

В дальнейшем полагаем $F=K^{*}\in\KK_f$. Вариационным методом
аналогично~\cite{PeRa94} доказывается, что экстремальный компакт $F$ (не
разбивающий плоскость) состоит из конечного числа кусочно-аналитических дуг,
является совокупностью критических траекторий некоторого квадратичного
дифференциала и обладает $S$-свойством~\eqref{spro2}.
Ниже (см.~лемму~\ref{l2}) приводится схема доказательство
этого результата
для случая {\it произвольного\/} параметра
$\theta\geq0$ при условии, что соответствующий экстремальный компакт
$F(\theta)$ существует, является допустимым, $F(\theta)\in\KK_f$, и
$F\not\ni\infty$; как показано выше в рассматриваемом здесь случае $\theta=1$
это условие выполняется. Подробное доказательство этого результата будет дано в
работе~\cite{RaSu11}.

Ниже мы в основном следуем схеме работы~\cite{PeRa94}.

\begin{lemma}\label{l2}
Пусть $F\in\KK_f$ -- стационарный компакт для функционала
энергии~\eqref{2.4} такой, что $F\not\ni\infty$. Тогда $F$ состоит из конечного числа кусочно аналитических
дуг, обладает $S$-свойством~\eqref{spro2} и
является замыканием критических траекторий некоторого квадратичного дифференциала.
\end{lemma}
\begin{proof}
Пусть $F\in\KK_f$ -- стационарный компакт для функционала энергии~\eqref{2.4}:
\begin{equation}
J(F,\lambda_F)=\max_{K\in\KK_f}J(K,\lambda_K),
\qquad\lambda_K\in M_1(E),
\label{v1}
\end{equation}
$\{a_1,\dots,a_{m-2}\}=F\cap{\Sigma_f}$ (точки $a_1,\dots,a_{m-2}\in\Sigma_f$ -- особые
точки функции $f$ многозначного характера, лежащие на $F$). Положим
$A_m(z)=\prod\limits_{j=1}^{m-2}(z-a_j)\cdot (z^2-1)$. Зафиксируем точку
$w\in\Omega=\myo{\CC}\setminus(E\cup{F})$ и рассмотрим следующее преобразование:
\begin{equation}
z_\tau=z+\tau\hh(z),
\quad\text{ где }
\hh(z)=\frac{A_m(z)}{(z-w)^{m+1}},
\label{v2}
\end{equation}
$\tau$ -- комплексный параметр. Фиксируем некоторую окрестность
$U_0=\{z:|z-w|<r_0\}$ точки $w$ такую, что $\myo U_0\subset\Omega$. Найдется
$\eps_0>0$ такое, что при всех $\tau$, $|\tau|<\eps_0$, отображение $z\mapsto
z_\tau$ взаимнооднозначно в области $\Omega_0=\Omega\setminus\myo{U}_0$ и
образ $F_\tau$ компакта $F$ -- допустимый компакт ($F_\tau\in\KK_f$).

Положим $D_\tau=\myo{\CC}\setminus{F_\tau}$.
Компакт $F$ состоит из конечного числа континуумов; следовательно, при
$|\tau|<\eps_0$ компакт $F_\tau$ также состоит из конечного числа континуумов.
Поэтому существуют функции Грина $g_{D}(z,\zeta)$ и $g_{D_\tau}(z,\zeta)$.
Найдем формулу для вариации функции Грина
\begin{equation}
\delta g(z,\zeta):=g_{D_\tau}(z,\zeta)-g_{D}(z,\zeta)
\label{v3}
\end{equation}
при $|\tau|<\eps_0$, $\eps_0$ -- достаточно мало, $z,\zeta$ -- фиксированы,
$z,\zeta\in{U}\Subset{D}$,
$D=\myo{\CC}\setminus{F}$.

Рассмотрим линию уровня $\Gamma$ функции Грина $g_{D}(t,z_0)$, где
$z_0\in{U}$ -- фиксированно, а $\Gamma$ ``достаточно близка'' к $F$. Ясно, что
при малых $\tau$ компакт $F_\tau$ лежит внутри этой линии уровня.
Воспользуемся вариационной формулой Адамара
(см.~\cite[приложение, \S~3, формула~\thetag{3.3}]{CoSh53}):
\begin{equation}
\delta g(z,\zeta)=-\frac1{2\pi}\int_\gamma
\frac{\partial g(t,z)}{\partial n_t}
\frac{\partial g(t,\zeta)}{\partial n_t}\delta n_t\,ds_t+O(\eps^2),
\label{v4}
\end{equation}
считая, что мы применяем ее для $\gamma=\Gamma$ и функции Грина для внешности
$\Gamma$ ($\delta n_t$ -- $\eps$-вариация вдоль нормали, в дальнейшем в
полученных ниже формулах надо будет сделать предельный переход при
$\Gamma\to{F}$).

Поскольку на линии уровня
$$
\frac{\partial g(t,z)}{\partial n_t}=G'(t,z)e^{i\alpha_t},
$$
где $G(t,z)=g(t,z)+ig^{*}(t,z)$ -- комплексная функция Грина и производная
берется по первому аргументу, то для вариации~\eqref{v1} из формулы
Адамара~\eqref{v4} получаем:
\begin{align}
\delta g(z,\zeta)
&=\Re\biggl\{\frac\tau{2\pi i}
\int_\gamma G'(t,z)G'(t,\zeta)\hh(t)\,dt\biggr\}+O(\eps^2)\notag\\
&=\Re\biggl\{\frac\tau{2\pi i}
\int_{|\zeta-w|=r_0}G'(t,z)G'(t,\zeta)\frac{A_m(t)}{(t-w)^{m+1}}\,dt\biggr\}
\notag\\
&\qquad\Re\biggl\{\frac\tau{2\pi i}
\int_{|t-z|=r_0}G'(t,z)G'(t,\zeta)\frac{A_m(t)}{(t-w)^{m+1}}\,dt\biggr\}\notag\\
&\qquad\Re\biggl\{\frac\tau{2\pi i}
\int_{|t-\zeta|=r_0}G'(t,z)G'(t,\zeta)\frac{A_m(t)}{(t-w)^{m+1}}\,dt\biggr\}
+O(\eps^2)\notag
\end{align}
Следовательно,
\begin{align}
\delta g(z,\zeta)
&=\Re\biggl\{\tau\Res_{t=w}\biggl(
G'(t,z)G'(t,\zeta)\frac{A_m(t)}{(t-w)^{m+1}}\biggr)\biggr\}\notag\\
&\qquad\Re\biggl\{\tau\Res_{t=z}\biggl(
G'(t,z)G'(t,\zeta)\frac{A_m(t)}{(t-w)^{m+1}}\biggr)\biggr\}\notag\\
&\qquad\Re\biggl\{\tau\Res_{t=\zeta}\biggl(
G'(t,z)G'(t,\zeta)\frac{A_m(t)}{(t-w)^{m+1}}\biggr)\biggr\}
+O(\eps^2)\notag\\
&=\Re\biggl\{\tau\frac1{m!}
\biggl(G'(w,z)G'(w,\zeta)A_m(w)\biggr)^{(m)}_{w}\notag\\
&\qquad
-\tau\frac{G'(z,\zeta)A_m(z)}{(z-w)^{m+1}}
-\tau\frac{G'(\zeta,z)A_m(\zeta)}{(\zeta-w)^{m+1}}
\biggr\}+O(\eps^2).
\label{v5}
\end{align}
Из~\eqref{v5} получаем, что
\begin{align}
\delta_\tau g(z,\zeta):
&=g_{D_\tau}(z_\tau,\zeta_\tau)-g_{D}(z,\zeta)\notag\\
&=\Re\biggl\{\tau\frac1{m!}
\biggl(G'(w,z)G'(w,\zeta)A_m(w)\biggr)^{(m)}_{w}\biggr\}+O(\eps^2).
\label{v5.2}
\end{align}
Аналогично
\begin{align}
\delta_\tau\log\frac1{|z-\zeta|}:
&=\log\frac1{|z_\tau-\zeta_\tau|}-\log\frac1{|z-\zeta|}\notag\\
&=\Re\biggl\{\tau\frac1{m!}
\biggl(\frac{A_m(w)}{(z-w)(\zeta-w)}\biggr)^{(m)}_{w}\biggr\}+O(\eps^2).
\label{v5.3}
\end{align}
Вариации~\eqref{v5.3} логарифмического ядра и~\eqref{v5.2} функции Грина
порождают следующую вариацию функционала
энергии~\eqref{2.4} для стационарного компакта
\begin{equation}
\delta_\tau J(F,\lambda_F)=\eps\Re\biggl\{e^{i\pfi}
\frac1{m!}
\biggl(A_m(w)\Bigl[\theta\bigl(\myh{\lambda}(w)\bigr)^2
+\bigl(\sG'_{\lambda,F}(w)\bigr)^2\Bigr]\biggr)^{(m)}_w\biggr\}
+O(\eps^2),
\label{v6}
\end{equation}
где
$$
\myh{\lambda}(w)=\sV'_\lambda(w)
=\int_{E}\frac{d\lambda(x)}{w-x}
,\qquad
\sG_{\lambda,F}(w)=\int_{E}G(w,x)\,d\lambda(x),
$$
$\tau=\eps e^{i\pfi}$ -- комплексный параметр. В силу
стационарности~\eqref{v1} компакта $F$ имеем:
$\delta_\tau J(F,\lambda_F)\leq0$ при
всех достаточно малых $\tau$. Поэтому из~\eqref{v6} получаем
\begin{equation}
\biggl(A_m(w)\bigl[\theta\bigl(\myh{\lambda}(w)\bigr)^2
+\bigl(\sG'_{\lambda,F}(w)\bigr)^2\bigr]\biggr)^{(m)}_w
\equiv0.
\label{v7}
\end{equation}
Следовательно,
\begin{equation}
A_m(w)\bigl[\theta\bigl(\myh{\lambda}(w)\bigr)^2
+\bigl(\sG'_{\lambda,F}(w)\bigr)^2\bigr]
\equiv B_{m-1}(w),
\label{v8}
\end{equation}
где $B_{m-1}$ -- полином степени $\leq{m-1}$. Непосредственно из~\eqref{v8}
вытекает следующее соотношение для стационарного компакта $F\in\KK_f$ и
соответствующей равновесной меры $\lambda=\lambda_F$ (с носителем на отрезке $[-1,1]$):
\begin{equation}
\label{v8.2}
\theta\bigl(\myh{\lambda}(w)\bigr)^2
+\bigl(\sG'_{\lambda,F}(w)\bigr)^2\equiv\frac{B_{m-1}(w)}{A_m(w)},
\end{equation}
а для гринова потенциала стационарного компакта $F$ имеем:
\begin{equation}
G^\lambda_F(z)=\int_{-1}^1 g_F(z,x)\,d\lambda(x)
=\Re\int_{a_1}^z\sqrt{\frac{B_{m-1}(\zeta)}{A_m(\zeta)}
-\theta\bigl(\myh{\lambda}(\zeta)\bigr)^2\,}\,d\zeta
\label{v9}
\end{equation}
(путь интегрирования в ~\eqref{v9} не пересекает отрезок $E$).
Из~\eqref{v9} вытекает справедливость утверждений леммы, в том числе --
$S$-свойство компакта $F$:
$$
\frac{\partial G^{\lambda}_{F}}{\partial n_{+}}(z)
=\frac{\partial G^{\lambda}_{F}}{\partial n_{-}}(z),
\qquad z\in F_0.
$$
Лемма~\ref{l2} доказана.
\end{proof}

\begin{remark}\label{r1}
1) Отметим, что утверждения леммы~\ref{l2} справедливы при любом значении параметра
$\theta\geq0$ (в формуле для функционала энергии $J(\mu;\theta)=\int(\theta
V^\mu(z)+G_F^\mu(z))\,d\mu(z)$)
при условии, что стационарный компакт $F=F(\theta)\in\KK_f$
{\it существует}.

2) Поскольку $F\not\ni\infty$, то непосредственно из~\eqref{v8}
вытекает, что $B_{m-1}=B_{m-2}$, где $\mdeg{B_{m-2}}\leq{m-2}$.

3) Наконец, при $\theta=0$ из~\eqref{v9}
получаем
\begin{equation}
G^\lambda_F(z)=\int_{-1}^1 g_F(z,x)\,d\lambda(x)
=\Re\int_{a_1}^z\sqrt{\frac{B_{m-2}(\zeta)}{A_m(\zeta)}\,}\,d\zeta.
\label{v10}
\end{equation}
\end{remark}

\subsection{}\label{s3s4}
Приведем теперь характеризацию компакта $F=F(1)$ с помощью двулистной римановой
поверхности, при этом окажется, что соответствующая постоянная равновесия
$w=w_F$ равна постоянной Робена для компакта $F^{(1)}$.

Пусть $K\in\KK_f$ -- произвольный допустимый компакт. Определим риманову
поверхность $\RS$ уравнением $w^2=z^2-1$, положим
$\RS=\RS^{(1)}\sqcup\RS^{(2)}\sqcup\Gamma$, где $\Gamma$~-- замкнутый цикл
на $\RS$, проходящий через точки $\pm1$ и соответствующий отрезку $[-1,1]$
при каноническом проектировании,
$\RS^{(1)}$ -- первый (открытый) лист римановй поверхности, $\RS^{(2)}$ --
второй лист.
Пусть функция $u(z^{(1)})=G^{\lambda}_K(z)$, $z^{(1)}\in\RS^{(1)}$. Тогда
из условий равновесия и симметрии относительно вещественной оси вытекает,
что функция $u$ гармонически продолжается на второй лист $\RS^{(2)}$ римановй
поверхности формулой $u(z^{(2)})=w-V^\lambda(z)$. Полученная функция $u(\zz)$,
$\zz\in\RS$, обладает следующими свойствами: $u$ -- гармоническая функция
на $\RS\setminus(K^{(1)}\cup\{\infty^{(2)}\})$, в точке $\zz=\infty^{(2)}$
функция $u$ имеет логарифмическую особенность, $u(\zz)=0$ при $\zz\in K^{(1)}$.
Следовательно, $u(\zz)=g_{K^{(1)}}(\zz,\infty^{(2)})$ -- функция Грина области
$D=\RS\setminus{K^{(1)}}$ с особенностью в бесконечно удаленной точке
$\zz=\infty^{(2)}$, $w_K=\gamma^{(2)}$ -- соответствующая постоянная Робена,
$c(K)=e^{-w_K}$ -- емкость множества $K^{(1)}$ (на $\RS$ относительно точки
$\zz=\infty^{(2)}$). Тем самым задача~\eqref{51.8} о максимуме постоянной
равновесия $w_K$ в классе $K\in\KK_f$ эквивалентна задаче о минимуме емкости
$c(K)$ в этом классе. Соответствующий экстремальный компакт $\myt{F}^{(1)}$
является компактом минимальной емкости (на $\RS$ относительно точки
$\zz=\infty^{(2)}$),
определяется однозначно и вполне характеризуется $S$-свойством:
\begin{equation}
\frac{\partial g_{\myt F^{(1)}}(\zeta^{(1)},\infty^{(2)})}{\partial n_+}=
\frac{\partial g_{\myt F^{(1)}}(\zeta^{(1)},\infty^{(2)})}{\partial n_-},
\qquad\zeta^{(1)}\in\myt{F}^{(1)}_0
\label{8.1}
\end{equation}
(см.~\cite{Sta85a}--\cite{Sta85c},~\cite{GoRa87}). В силу определения функции
$u$, соотношение~\eqref{8.1} эквивалентно $S$-свойству~\eqref{spro2}. Тем
самым, $F=\myt{F}$ и $G^\lambda_F(z)=g_{\Phi(F)}(\Phi(z),0)$ при
$z\in\myo{\CC}\setminus{E}$, т.е. стационарный в смысле экстремальной
задачи~\eqref{51.9} компакт $F$ совпадает с образом компакта Шталя.

\section{Доказательство теоремы~\ref{t2}}\label{s4}

\subsection{}\label{s4s1}
Перейдем теперь непосредственно к изучению сходимости (нелинейных) АПЧ.

Приведем сначала некоторые необходимые нам свойства общих полиномов Фабера,
соответствующих произвольному континууму~$E\Subset\CC$ (см.~\cite{SmLe64}).

Пусть $E$~-- произвольный континуум в~$\CC$, не разбивающий плоскость и не
вырождающийся в точку, $\Omega=\myo{\CC}\setminus{E}$. Положим $w=\Phi(z)$ -- функция, конформно
(и однолистно) отображающая область $\Omega$ на внешность единичного круга
$\{w:|w|>1\}$ так, что $\Phi(\infty)=\infty$, $\Phi'(\infty)>0$,
$\Phi\colon\Omega\to\{w:|w|>1\}$, $z=\Psi(w)$~-- обратная функция. Для произвольного
$\rho>1$ положим $\Gamma_\rho=\{z:|\Phi(z)|=\rho\}$ -- линия уровня
отображающей функции $\Phi$, $D_\rho$ -- внутренность кривой $\Gamma_\rho$,
$D_1:=E$. Тогда полиномы Фабера $\Phi_n(z)=\Phi_n(z;E)$ определяются по формуле
\begin{equation}
\Phi_n(z):=\frac1{2\pi i}\int_{\Gamma_\rho}\frac{\Phi^n(\zeta)}{\zeta-z}d\zeta,
\qquad z\in D_\rho,\quad\rho>1
\label{9.1}
\end{equation}
($\Phi_n(z)$ -- главная часть разложения функции $\Phi^n(z)$ в ряд Лорана
в бесконечно удаленной точке $z=\infty$).

Из~\eqref{9.1} вытекает, что $\Phi_n(z)=\Phi^n(z)(1+o(1))$,
$n\to\infty$, локально равномерно в~$\Omega$. Следовательно, существует
\begin{equation}
\lim_{n\to\infty}|\Phi_n(z)|^{1/n}=|\Phi(z)|,\qquad z\in \Omega.
\label{9.2}
\end{equation}
Из~\eqref{9.1} и~\eqref{9.2} легко вытекает, что если функция $f$ --
голоморфная на $E$, $f\in\sH(E)$, то $f$ разлагается в ряд Фабера
(ср.~\eqref{1.01})
\begin{equation}
f(z)=\sum_{k=0}^\infty c_k\Phi_k(z),\qquad c_k=c_k(f),
\label{9.4}
\end{equation}
сходящийся к $f$ локально равномерно в канонической области
$D_0(f):=D_{\rho_0(f)}$, где
\begin{equation}
\frac1{\rho_0(f)}=\varlimsup_{k\to\infty}|c_k|^{1/k}<1;
\label{9.5}
\end{equation}
область $D_0(f)$ -- максимальная каноническая область голоморфности
функции~$f$.

В классе функций $\sH(E)$ стандартным образом определим (линейный) {\it
оператор Фабера} $\sU\colon\sH(E)\to\sH(\myo{\DD})$, $\DD:=\{w:|w|<1\}$,
по формуле
\begin{equation}
\sU(f)(w)=\sum_{k=0}^\infty c_k w^k,
\label{9.6}
\end{equation}
где $f\in\sH({E})$ задана разложением
$f(z)=\sum\limits_{k=0}^\infty c_k\Phi^k(z)$. В силу~\eqref{9.5} функция
$\myt{f}=\sU(f)\in\sH(\myo{\DD})$, при этом $\rho_0(f)=R_0(\myt f)$,
где $R_0$ -- радиус голоморфности соответствующей функции. Тем самым
$\sU\colon\sH(E)\to\sH(\myo{\DD})$. Нетрудно видеть, что отображение $\sU$
ограниченно, биективно и справедлива следующая формула обращения
\begin{equation}
f(z)=\sU^{-1}(\myt f)(z):=\frac1{2\pi i}\int_{\Gamma_\rho}
\frac{\myt f(\Phi(\zeta))}{\zeta-z}\,d\zeta,\qquad
z\in D_\rho,\quad 1<\rho<\rho_0(f).
\label{9.7}
\end{equation}

Для произвольных $n,m\in\NN_0$ положим $\sR_{n,m}(E)=\sR_{n,m}\cap\sH(E)$,
$\sR_{n,m}(\myo{\DD})=\sR_{n,m}\cap\sH(\myo{\DD})$. Для произвольных
фиксированных $m\in\NN$ и $\rho>1$ обозначим через
$\sM_m(D_\rho)$ класс функций $f\in\sH(E)$, допускающих мероморфное
продолжение в каноническую область $D_\rho$ и имеющих там ровно\footnote{Как
обычно нули и полюсы функций считаются с учетом их кратностей.}
$m$ полюсов; аналогичным образом определяется и
$\sM_m(K_\rho)$ для $\myt f\in\sH(\myo{\DD})$, где $K_\rho=\{w:|w|<\rho\}$.

Справедливо следующее утверждение (см.~\cite{Sue80}, а также~\cite{Ged81}
и~\cite{Sue09b}) о свойствах оператора Фабера.

\begin{lemma}\label{l3}
Оператор Фабера~$\sU$ отображает:

1) при $n\geq m-1$ множество $\sR_{n,m}(E)$ линейно и взаимно однозначно на множество
$\sR_{n,m}(\myo{\DD})$, причем точка $z_0$ является полюсом кратности
$\mu\geq1$ функции $r(z)\in\sR_{n,m}(E)$ тогда и только тогда, когда
точка $w_0=\Phi(z_0)$ является полюсом той же кратности $\mu$ функции
$R(w)=\sU(r)(w)\in\sR_{n,m}(\myo{\DD})$;

2) множество $\sM_m(D_\rho)$, $\rho>1$, линейно и взаимно однозначно на
множество $\sM_m(K_\rho)$, причем точка $z_0$ является полюсом кратности
$\mu\geq1$ функции $f(z)\in\sM_m(D_\rho)$ тогда и только тогда, когда точка
$w_0=\Phi(z_0)$ является полюсом той же кратности $\mu$ функции
$\myt f(w)=\sU(f)(w)$.
\end{lemma}

Эта простое утверждение впервые было сформулировано в~\cite{Sue80} (см.
также~\cite{Sue09b}). Позднее эти
свойства оператора Фабера и вытекающие из них следствия (в том числе, аналог
теоремы Монтессу) были переоткрыты другими
авторами~\cite{Ell83},~\cite{ElSa84}.

Вышеуказанные свойства оператора Фабера понадобятся нам здесь для случая,
когда континуум $E=[-1,1]$, $\Psi(w)=\Zh(w)=(w+1/w)/2$ -- функция Жуковского,
$w=\Phi(z)=z+\sqrt{z^2-1}$ и выбрана такая ветвь корня, что $\Phi(z)\sim2z$
при $z\to\infty$.

\subsection{}\label{s4s2}
Из следующего соотношения (см.~\eqref{2.03}, определение нелинейных АПЧ)
\begin{equation}
(f-F_n)(z)=cT_{2n+1}(z)+\dotsb
\label{9.8}
\end{equation}
вытекает, что
\begin{equation}
(\myt{f}-\myt{F}_n)(w)=cw^{2n+1}+\dotsb,
\label{9.9}
\end{equation}
где в силу леммы~\ref{l3} функция $\myt{F}_n\in\sR_n(\myo{\DD})$.
Следовательно, рациональная функция $\myt{F}_n(w)=U(F_n)(w)=[n/n]_{\myt{f}}(w)$
-- диагональная аппроксимация Паде ряда $\sum\limits_{k=0}^\infty c_kw^k$
порядка~$n$.

Нетрудно видеть, что оператор Фабера (в нашем случае -- Фабера--Чебышёва)
сохраняет характер особенностей (в частности -- ветвления) при переходе от
функций $f$ к функции $\myt{f}$ и устанавливает взаимно однозначное
соответствие между допустимыми компактами для функций $f$ и $\myt{f}$ (заданных
в плоскости $\CC_z$ и $\CC_w$ соответственно). Следовательно, при условии, что
$\myt{f}\in\HH(\infty)$, формулу~\eqref{9.7} можно записать в виде
\begin{equation}
f(z)=R(z)+\frac1{2\pi i}\oint_F \frac{[\myt{f}](\Phi(\zeta))}{\zeta-z}d\zeta,
\qquad z\in\Omega\setminus{E},
\label{9.10}
\end{equation}
где $F\in\KK_f$ -- $S$-симметричный компакт для $f$, $R=U^{-1}(r)$, $r$ --
сумма главных частей функции $\myt{f}$ в $\myo\CC_w\setminus{\Phi(F)}$,
$[\myt{f}]$ -- голоморфная составляющая функции $\myt{f}$ в
$\myo\CC_w\setminus{\Phi(F)}$, и под $\displaystyle\oint_{F}$ понимается
интеграл по любому контуру, охватывающему $F$ и отделяющему $F$ от точки $z$ и
отрезка~$E$. Аналогичное~\eqref{9.10} представление справедливо и для
диагональных АПЧ $F_n$ функции $f$.

Из первой теоремы Шталя (см.~\cite{Sta86a}, а также приложение~\ref{apl})
вытекает, что $[n/n]_{\myt{f}}\overset\mcap\rightarrow{\myt{f}}$ на компактных
подмножествах области $\myo{\CC}_w\setminus\myt{F}$, где $\myt{F}=\Phi(F)$
-- компакт Шталя\footnote{Более точно, компакт Шталя есть множество
$S=\{z:1/z\in\myt{F}\}$ и определяется относительно бесконечно удаленной
точки.}
для функции $\myt{f}$. При этом
\begin{equation}
|\myt{f}(w)-[n/n]_{\myt{f}}(w)|^{1/n}\overset\mcap\rightarrow
e^{-2g_{\myt{F}}(w,0)},\qquad w\in{\myt{K}}\Subset{\myt{D}}\setminus{E},
\label{9.11}
\end{equation}
$\myt{D}=\myo\CC_w\setminus{\myt{F}}$, и для $\myt{Q}_n$ -- знаменателей АП
функции $\myt{f}$ имеем:
\begin{equation}
\frac1n\mu(\myt{Q}_n)\to\lambda_{\myt{F}}\quad\text{ -- равновесная мера
для $\myt{F}$}.
\label{9.12}
\end{equation}
Из приведенной леммы~\ref{l3} и соотношений~\eqref{9.8}--\eqref{9.12} вытекает,
что $F_n\overset\mcap\rightarrow f$ на компактных
подмножествах области $\myo{\CC}_z\setminus{F}$, при этом в силу~\eqref{9.11}
и равенства $G_F^\lambda(z)=g_{\myt{F}}(\Phi(z),0)$,
$z\in\myo{\CC}\setminus{E}$, имеем
\begin{equation}
|f(z)-F_n(z)|^{1/n}\overset\mcap\rightarrow
e^{-2G^{\lambda}_{F}(z)},\qquad z\in{K}\Subset{D},
\label{9.13}
\end{equation}
и все нули полиномов $Q_n$ за исключением $o(n)$ из них притягиваются к
компакту~$F$. Сходимость диагональных АПЧ функции~$f$ доказана.
Отметим, что из сказанного выше вытекает единственность диагональных АПЧ
(по соответствующей подпоследовательности) для достаточно больших~$n$.

\subsection{}\label{s4s3}
Докажем теперь утверждения о предельном распределении нулей $Q_n$ и точек
интерполяции.

Итак, пусть нелинейная АПЧ $F_n$ существует при всех $n\geq n_0$, принадлежащих
некоторой бесконечной последовательности $\Lambda\subset\NN$, и $F_n\in\HH(E)$.
Все дальнейшие рассуждения будут проводится для таких~$n$.
Из~\eqref{2.02} получаем следующие соотношения ортогональности:
\begin{equation}
\int_{E}(f-F_n)(x)T_k(x)\,d\tau(x)=0,\qquad k=0,1,\dots,2n.
\label{9.14}
\end{equation}
Функции $f$ и $F_n$ вещественнозначны на $E$. Поэтому из соотношений~\eqref{9.14}
вытекает, что разность $f-F_n$ обращается в нуль на $E$ по-крайней мере в
$(2n+1)$-й точке. Пусть $\omega_{2n}(z)=z^{2n+1}+\dotsb$ -- полином с нулями в
этих точках. Положим $F_n=P_n/Q_n$, тогда имеем: функция $(Q_nf-P_n)/\omega_{2n}$
голоморфна на $E$. Пусть $s$ -- полином фиксированной степени $m$ с единичным
старшим коэффициентов, часть нулей которого совпадает с полюсами функции~$f$ в
$D=\myo{\CC}\setminus{F}$ и имеет кратность, равную кратности соответствующего
полюса. Кроме того, будем
считать, что среди нулей $s$ содержатся все точки ветвления функции $f$
с такой кратностью, что функция $s\cdot\Delta{f}=s\cdot(f_{+}-f_{-})$ --
непрерывна на $F$; здесь и далее через $\Delta{f}=f_{+}-f_{-}$ обозначается
скачок функции
$f$ на дугах, составляющих компакт $F$. Ясно, что $s$ -- вещественный полином,
не имеющий нулей на $E$.
Функция $s(Q_nf-P_n)/\omega_{2n}$ голоморфна в $D$. Теперь стандартным способом
(см., например,~\cite{GoLo78},~\cite{GoRa87}) получаем следующие соотношения:
\begin{spreadlines}{6mm}
\begin{gather}
(f-F_n)(z)=\frac{\omega_{2n}(z)}{Q_n(z)q(z)s(z)}\frac1{2\pi i}
\oint_{F}\frac{Q_n(t)q(t)s(t)f(t)\,dt}{\omega_{2n}(t)(t-z)},\qquad z\in D,
\label{9.15}\\
\oint_{F}\frac{Q_n(t)t^ks(t)f(t)\,dt}{\omega_{2n}(t)}=0,\qquad k=0,1,\dots,n-1-m,
\quad m=\mdeg{s},\label{9.16}\\
\int_{-1}^{1}\omega_{2n}(x)x^j\frac1{Q_n(x)q(x)s(x)}
\biggl\{\oint_{F}\frac{Q_n(t)q(t)s(t)f(t)\,dt}{\omega_{2n}(t)(t-x)}\biggr\}
d\tau(x)=0,\quad j=\overline{0,2n},
\label{9.17}
\end{gather}
\end{spreadlines}
где как и выше под $\displaystyle\oint_{F}$ понимается интеграл
по любому контуру, охватывающему $F$ и отделяющему $F$ от точки $z$ и
отрезка~$E$, $q$ -- произвольный полином степени $\leq{n-m}$. Так как
функция $s\cdot\Delta{f}$ непрерывна на $F$, то в
соотношениях~\eqref{9.15}--\eqref{9.17} от интеграла $\displaystyle\oint_{F}$
можно перейти к интегралу $\displaystyle\int_{F}$, а $f(t)$ заменить на
$\Delta{f}(t)$. Кроме того, вместо полинома $q$ в эти соотношения можно
подставить полином $\myt{Q}_n$, который отличается от $Q_n$
сомножителем, степень которого $\leq{m}$. Таким образом, получаем
\begin{spreadlines}{6mm}
\begin{gather}
(f-F_n)(z)=\frac{\omega_{2n}(z)}{Q_n(z)\myt{Q}_n(z)s(z)}\frac1{2\pi i}
\int_{F}\frac{Q_n(t)\myt{Q}_n(t)s(t)\Delta f(t)\,dt}{\omega_{2n}(t)(t-z)},\qquad z\in D,
\label{9.18}\\
\int_{F}\frac{Q_n(t)t^ks(t)\Delta f(t)\,dt}{\omega_{2n}(t)}=0,\qquad
k=0,1,\dots,n-1-m, \label{9.19}\\
\int_{-1}^{1}\omega_{2n}(x)x^j\frac1{Q_n(x)\myt{Q}_n(x)s(x)}
\biggl\{\int_{F}\frac{Q_n(t)\myt{Q}_n(t)s(t)\Delta f(t)\,dt}{\omega_{2n}(t)(t-x)}\biggr\}
d\tau(x)=0,\quad j=\overline{0,2n}.
\label{9.20}
\end{gather}
\end{spreadlines}
В силу \eqref{9.13} справедливо соотношение
\begin{equation}
\label{9.201}
|f(z)-F_n(z)|^{1/2n}\overset\mcap\rightarrow
e^{-G^{\lambda}_{F}(z)}.
\end{equation}
Так как все нули полинома $\omega_{2n}$ лежат на отрезке $E$, то (переходя при
необходимости к подпоследовательности) имеем:
\begin{equation}
|\omega_{2n}(z)|^{1/2n}\to e^{-V^{\mu_\omega}(z)}
\label{9.21}
\end{equation}
локально равномерно в $\CC\setminus{E}$, где $\mu_\omega\in M_1(E)$. Поскольку
почти все нули полиномов $Q_n$ и $\myt{Q}_n$ (кроме $o(n)$ из них)
притягиваются к компакту $F$, то
\begin{equation}
|Q_n(z)|^{1/n}\overset\mcap\rightarrow e^{-V^{\mu_Q}(z)},
\quad
|\myt Q_n(z)|^{1/n}\overset\mcap\rightarrow e^{-V^{\mu_Q}(z)},
\qquad
z\in\myo{\CC}\setminus{F},
\label{9.22}
\end{equation}
где $\mu_Q\in M_1(F)$. Следовательно, из~\eqref{9.201}--\eqref{9.22} получаем:
\begin{equation}
V^{\mu_Q}(z)-V^{\mu_\omega}(z)+V(z)=G^\lambda_F(z),\qquad
z\in\myo{\CC}\setminus(E\cup F),
\label{9.23}
\end{equation}
где
\begin{equation}
V(z):=\lim_{n\to\infty}\log
\biggl|\int_{F}\frac{Q_n(t)\myt{Q}_n(t)s(t)\Delta{f}(t)\,dt}{\omega_{2n}(t)(t-z)}
\biggr|^{1/2n}.
\label{9.24}
\end{equation}
Функция $V$ -- гармоническая вне $E\cup{F}$ и субгармоническая на $E$.
Из соотношения равновесия~\eqref{9.23} вытекает, что разность
$V-V^{\mu_\omega}$ непрерывно продолжается на $E$. Потенциал $V^{\mu_\omega}$ -- супергармоническая функция
на $E$, тем самым $V^{\mu_\omega}$ полунепрерывна снизу на $E$. Функция $V$
полунепрерывна сверху на $E$ и отличается от $V^{\mu_\omega}$ на непрерывную
функцию. Следовательно, обе функции $V$ и $V^{\mu_\omega}$ непрерывны на
$E$. Теперь уже из соотношения~\eqref{9.20}, используя вещественность всех
входящих в это соотношение функций, стандартным методом работы~\cite{GoRa84}
получаем следующее соотношение равновесия на $E$:
\begin{equation}
-2V^{\mu_\omega}(x)+V^{\mu_Q}(x)+V(x)\equiv\const,\qquad x\in{E}.
\label{9.25}
\end{equation}
Из соотношений~\eqref{9.23} и~\eqref{9.25} вытекает, что
\begin{equation}
V^{\mu_\omega}(x)+G^\lambda_{F}(x)\equiv\const,\qquad x\in{E}.
\label{9.26}
\end{equation}
Следовательно, в силу условий равновесия~\eqref{51.1} имеем:
$$
V^{\mu_\omega}(x)\equiv V^\lambda(x)+\const,\qquad x\in{E},
\quad \mu_\omega,\lambda\in M_1(E).
$$
Отсюда получаем, что $\mu_\omega=\lambda$.

Наконец, воспользуемся соотношениями ортогональности~\eqref{9.19}. Из
$S$-свойства компакта $F$ и равенства $\mu_\omega=\lambda$ вытекает, что
мы находимся в условиях общей теоремы~3 работы~\cite{GoRa87}. Непосредственно
из этой теоремы для рассматриваемого здесь частного случая
$\psi(z)=-V^\lambda(z)$ вытекает, что $\mu_Q=\myt{\lambda}$ -- выметание меры
$\lambda$ на $F$, а $V(z)\equiv\const$.

Из сказанного выше вытекает, что соотношения
\begin{equation}
\frac1{2n}\mu(\omega_{2n})\to\lambda,\qquad
\frac1{n}\mu(Q_n)\to\myt\lambda
\label{9.27}
\end{equation}
имеют место для любой подпоследовательности $\Lambda\subset\NN$. Следовательно,
$\lambda$ и $\myt\lambda$ -- единственные предельные точки для
последовательностей нормированных мер $\dfrac1{2n}\mu(\omega_{2n})$ и $\dfrac1{n}\mu(Q_n)$
соответственно. Тем самым все утверждения теоремы~2 доказаны.

\section{Приложение: теоремы Шталя}\label{apl}

\subsection{Первая теорема Шталя}\label{as1}
Пусть функция $f\in\HH(\infty)$ задана в точке $z=\infty$ сходящимся рядом
Лорана
\begin{equation}
f(z)=\sum_{k=0}^\infty\frac{c_k}{z^{k+1}}.
\label{alor1}
\end{equation}
Предположим, что существует конечное непустое множество точек
$\Sigma\subset\CC$ такое, что функция $f$ продолжается
из окрестности бесконечно удаленной точки по любому
пути, лежащему в $\CC\setminus{\Sigma}$, и
$f\notin\MM(\myo{\CC}\setminus{\Sigma})$. В таком случае будем писать
$f\in\myA_\infty(\myo{\CC}\setminus{\Sigma})$.

Обозначим через $[n/n]_f$ диагональную аппроксимацию Паде функции
$f$ (в точке $z=\infty$): $[n/n]_f=P_n/Q_n$, где
$\mdeg{P_n},\mdeg{Q_n}\leq{n}$, $Q_n\not\equiv0$ и выполняется соотношение
\begin{equation}
(Q_nf-P_n)(z)=O\(\frac1{z^{n+1}}\),\quad z\to\infty.
\label{asta1-01}
\end{equation}

Для произвольного компакта $K\subset\myo{\CC}$ через $D_K$ обозначим связную
компоненту дополнения к $K$, содержащую точку $z=\infty$.
Пусть $\KK_f$ -- семейство компактов $K\subset\myo{\CC}$ таких, что
область $D_K$ регулярна относительно решения задачи Дирихле и $f$ допускает
мероморфное (однозначное аналитическое) продолжение из окрестности точки
$z=\infty$ в область $D_K$: $f\in\MM(D_K)$. Компакты $K\in\KK_f$ будем называть
допустимыми для функции $f\in\myA_\infty(\myo{\CC}\setminus{\Sigma})$.

Через $\mu(Q)$ обозначим меру, ассоциированную с произвольным полиномом~$Q$:
$$
\mu(Q)=\sum_{\zeta:Q(\zeta)=0}\delta_\zeta,
$$
где $\delta_\zeta$ -- мера Дирака с носителем в точке $\zeta$.

Следующий результат о сходимости диагональных аппроксимаций Паде принадлежит
Г.~Шталю~\cite{Sta85a}--\cite{Sta86b}.

\begin{theoremFsSt}
Пусть $f\in\myA_\infty(\CC\setminus\Sigma)$. Тогда

1) существует единственный компакт $F=F(f)$ такой, что $F\in\KK_f$~и
\begin{equation}
\mcap{F}=\min_{K\in\KK_f}\mcap{K},
\label{mincap}
\end{equation}
компакт $F$ состоит из конечного числа кусочно-аналитических дуг
и не разбивает плоскость;

2) для нулей $Q_n$ (знаменателей рациональных функций $[n/n]_f$) справедливо
предельное соотношение:
\begin{equation}
\frac1n\mu(Q_n)\to\lambda_F,\qquad n\to\infty,
\label{ameco1}
\end{equation}
где $\lambda_F$ -- (единичная) равновесная мера компакта $F$, сходимость мер
понимается в слабой топологии;

3) последовательность $[n/n]_f$ сходится по емкости внутри (на компактных
подмножествах) области $D=\myo\CC\setminus{F}$ к функции $f$:
\begin{equation}
\bigl|f(z)-[n/n]_f(z)\bigr|^{1/n}
\overset{\mcap}\longrightarrow e^{-2g_F(z)}<1,
\qquad
z\in{K},\quad K\subset{D},
\label{asta1-1}
\end{equation}
$g_F(z)=g_F(z,\infty)$ -- функция Грина для области $D$ с особенностью в бесконечно
удаленной точке.
\end{theoremFsSt}

Отметим, что область $D$ является областью максимальной сходимости
аппроксимаций Паде функции $f$;
скачок $\Delta f(\zeta)=(f^{+}-f^{-})(\zeta)$, $\zeta\in\ell$,
функции $f$ на любой открытой дуге $\ell\subset F$ отличен от тождественного
нуля: $\Delta f(\zeta)\not\equiv0$, $\zeta\in\ell$.

Соотношение~\eqref{asta1-1} означает следующее: {\it для любого компакта $K\Subset
D_f$ и любого числа $\varepsilon>0$ найдутся номер $n_0=n_0(K,\varepsilon)$
и число $\delta=\delta(K,\varepsilon)>0$ такие, что при всех $n\geq{n_0}$
справедливо соотношение:}
\begin{equation}
(e^{-2g_F(z)-\varepsilon})^{n}\leq|(f-[n/n]_f)(z)|
\leq(e^{-2g_F(z)+\varepsilon})^{n},
\qquad z\in K\setminus{e_n},\quad\mcap{e_n}<\delta.
\label{asta1-2}
\end{equation}

Из~\eqref{asta1-2} уже легко следует, что {\it каждый полюс $a$ функции $f$ в
области $D$ кратности $\nu(a)\ge1$ притягивает при $n\to\infty$ по-крайней мере
$\nu(a)$ полюсов рациональной функции $[n/n]_f$.}

Метод Шталя основан на том, что экстремальный компакт $F$
(см.~\eqref{mincap}) является совокупностью замыканий критических траекторий
некоторого квадратичного дифференциала, состоит из конечного числа
кусочно-аналитических дуг, не разбивает плоскость и вполне характеризуется
следующим свойством $S$-симметрии (или $S$-свойством):
\begin{equation}
\frac{\partial g^{}_F(\zeta,\infty)}{\partial n_{+}}=
\frac{\partial g^{}_F(\zeta,\infty)}{\partial n_{-}},
\qquad \zeta\in F_0,
\label{aspro1}
\end{equation}
$F_0$ -- совокупность открытых аналитических дуг, составляющих~$F$,
$\partial n_{\pm}$ -- производные по нормали с противоположных сторон~$F_0$.
Свойство неограниченной продолжаемости функции~$f$ в $\myo\CC\setminus{\Sigma}$
не существенно; важен лишь тот факт, что $f$~имеет ``правильный'' скачок
на некотором компакте $F$, удовлетворяющем условию~\eqref{aspro1}.

Компакты, обладающие тем или иным $S$-свойством, часто возникают в различных
задачах теории аппроксимаций и геометрической теории функций.

В частности, свойством симметрии обладают экстремальные кривые целого ряда
классических экстремальных задач геометрической теории функций; в этой теории
экстремальные кривые обычно описываются как траектории некоторого квадратичного
дифференциала
(см.~\cite{Jen62},~\cite{Gol66},~\cite{Kuz80}, \cite{Dub94a},~\cite{Dub94b}).

Первая теорема Шталя доказывается по следующей схеме.
Вначале на основе достаточно простых геометрических соображений (ср. лемма~\ref{l1})
устанавливается, что существует допустимый компакт $F\in\KK_f$ такой, что
\begin{equation}
\mcap{F}=\min_{K\in\KK_f}\mcap{K}.
\label{asta1-3}
\end{equation}
Затем с помощью вариационного метода доказывается, что $F$ является
замыканием критических траекторий некоторого квадратичного дифференциала.
Отсюда уже вытекает $S$-свойство~\eqref{aspro1}. Непосредственно на основе
этого $S$-свойства устанавливается предельное соотношение~\eqref{ameco1} для
нормированных знаменателей диагональных аппроксимаций Паде $[n/n]_f$.
Соотношение~\eqref{ameco1} влечет сходимость по емкости последовательности
$\{[n/n]_f\}$ к функции~$f$ внутри (на компактных подмножествах) области~$D$.
Единственность компакта $F$, удовлетворяющего условию~\eqref{asta1-3},
вытекает из сходимости~\eqref{ameco1}.

Соотношение~\eqref{ameco1} эквивалентно соотношению
\begin{equation}
\bigl|Q_n(z)\bigr|^{1/n}
\overset{\mcap}\longrightarrow C e^{g_{F}(z,\infty)},\qquad z\in D,\quad
C=\mcap F,
\label{apoco1}
\end{equation}
где старший коэффициент полинома $Q_n$ равен единице.

\subsection{Вторая теорема Шталя}\label{as2}
Пусть $E$ -- произвольный односвязный континуум в комплексной плоскости $\CC$,
функция $f$ голоморфна на $E$, $f\in\HH(E)$.
Предположим, что существует непустое конечное множество точек
$\Sigma\subset\myo{\CC}$ такое, что функция $f$
продолжается с компакта $E$ по любому пути, не пересекающему множество
$\Sigma$, и $f\not\in\MM(\myo\CC\setminus\Sigma)$;
в таком случае будем писать $f\in\myA_{E}(\myo\CC\setminus\Sigma)$.

Обозначим через $R_n=P_n/Q_n$ наилучшую равномерную рациональную аппроксимацию
функции~$f$ на $E$ в классе $\sR_n$ -- рациональных функций вида $r=p/q$,
$\mdeg{p},\mdeg{q}\leq{n}$, ${q}\not\equiv0$:
$$
\|f-R_n\|_{E}=\min_{r\in\sR_n}\|f-r\|_{E},
$$
$\|\cdot\|_{E}$ -- $\sup$-норма на $E$.

Для произвольного компакта $K\subset\myo{\CC}\setminus{E}$ через
$D_K$ обозначим связную компоненту дополнения к $K$, содержащую континуум~$E$.
Пусть $\KK_f$ -- семейство компактов $K$ таких, что
$K\subset\myo{\CC}\setminus{E}$, область $D_K$ регулярна относительно решения
задачи Дирихле и $f$ допускает мероморфное (однозначное аналитическое)
продолжение с континуума~$E$ в $D_K$: $f\in\MM(D_K)$. Компакты $K\in\KK_f$
будем называть допустимыми для функции $f\in\myA_{E}(\myo\CC\setminus\Sigma)$.

Для произвольного компакта $K\subset\myo{\CC}\setminus{E}$ такого, что область
$D_k$ регулярна относительно решения задачи Дирихле, и произвольной
положительной меры $\mu$ с носителем на $E$ положим
$$
G^\mu_K(z)=\int_{E}g_K(z,\zeta)\,d\mu(\zeta),\qquad
z\in\myo{\CC}\setminus(E\cup{K}),
$$
-- гринов (относительно $K$) потенциал меры $\mu$, $g_K(z,\zeta)$ -- функция
Грина для области $D_K$. Через $\myt\mu$ обозначим выметание меры $\mu$ из
области $D_K$ на $K$.

Через $\mu(Q)$ обозначим меру, ассоциированную с произвольным полиномом~$Q$:
$$
\mu(Q)=\sum_{\zeta:Q(\zeta)=0}\delta_\zeta,
$$
где $\delta_\zeta$ -- мера Дирака с носителем в точке $\zeta$.

Следующий результат о сходимости наилучших равномерных рациональных
аппроксимаций функции $f\in\myA_E(\myo{\CC}\setminus{\Sigma})$ принадлежит
Г.~Шталю~\cite{Sta85a}--\cite{Sta86b}.

\begin{theoremSeSt}
Пусть $f\in\myA_E(\myo\CC\setminus\Sigma)$. Тогда

1) существует единственный компакт $F=F(f)\in\KK_f$ такой, что
\begin{equation}
\mcap{(E,F)}=\min_{K\in\KK_f}\mcap{(E,K)},
\label{mincon}
\end{equation}
где $\mcap(E,K)$ -- емкость конденсатора $(E,K)$; компакт $F$ состоит из конечного числа кусочно-аналитических дуг
и не разбивает плоскость;

2) для нулей полинома $Q_n$ (знаменателя рациональной функции $R_n$)
справедливо предельное соотношение:
\begin{equation}
\frac1n\mu(Q_n)\to\myt\lambda,\qquad n\to\infty,
\label{ameco2}
\end{equation}
где $\lambda$ равновесная мера (с носителем на $E$) для гринова потенциала
компакта $F$, $G^\lambda_F(z)\equiv\const$ на $E$, $\myt\lambda$ -- выметание
меры $\lambda$ на $F$, сходимость мер понимается в слабой
топологии;

3) последовательность $R_n$ сходится по емкости внутри (на компактных
подмножествах) области $D=\myo{\CC}\setminus{F}$ и на компактных подмножествах
области $\Omega=\myo{\CC}\setminus({F}\cup{E})$ имеем:
\begin{equation}
\bigl|f(z)-R_n(z)\bigr|^{1/n}
\overset\mcap\longrightarrow
e^{-2G^\lambda_F(z)}<1,
\qquad z\in{K},\quad K\subset\Omega.
\label{asta2-1}
\end{equation}
\end{theoremSeSt}

Отметим, что область $D$ является максимальной областью сходимости
последовательности $R_n$;
скачок $\Delta f(\zeta)=(f^{+}-f^{-})(\zeta)$, $\zeta\in\ell$,
функции $f$ на любой открытой дуге $\ell\subset F$ отличен от тождественного
нуля: $\Delta f(\zeta)\not\equiv0$, $\zeta\in\ell$.

Соотношение~\eqref{asta2-1} означает следующее: {\it для любого компакта
$K\subset\Omega$ и любого числа $\varepsilon>0$ найдутся номер $n_0=n_0(K,\varepsilon)$
и число $\delta=\delta(K,\varepsilon)>0$ такие, что при всех $n\geq{n_0}$
справедливо соотношение:}
\begin{equation}
(e^{-2G^\lambda_F(z)-\varepsilon})^{n}
\leq|(f-R_n)(z)|\leq(e^{-2G^\lambda_F(z)+\varepsilon})^{n},
\qquad z\in K\setminus{e_n},\quad\mcap{e_n}<\delta.
\label{asta2-2}
\end{equation}

Из~\eqref{asta2-2} уже легко следует, что {\it каждый полюс $a$ функции $f$ в области
$D_f$ кратности $\nu(a)\ge1$ притягивает при $n\to\infty$ по-крайней мере
$\nu(a)$ полюсов рациональной функции $R_n$.}

Метод Шталя основан на том, что экстремальный компакт $F=F(f)$
(см.~\eqref{mincon}) является совокупностью замыканий критических траекторий
некоторого квадратичного дифференциала, состоит из конечного числа
кусочно-аналитических дуг, не разбивает плоскость и вполне характеризуется
следующим свойством $S$-симметрии (или $S$-свойством):
\begin{equation}
\frac{\partial G^{\lambda}_F(\zeta)}{\partial n_{+}}=
\frac{\partial G^{\lambda}_F(\zeta)}{\partial n_{-}},
\qquad \zeta\in F_0,
\label{aspro2}
\end{equation}
$F_0$ -- совокупность открытых аналитических дуг, составляющих~$F$,
$\partial n_{\pm}$ -- производные по нормали с противоположных сторон~$F_0$.
Свойство неограниченной продолжаемости функции~$f$ в
$\myo\CC\setminus{\Sigma}$ не существенно; важен лишь тот факт, что $f$~имеет
``правильный'' скачок на компакте~$F$, обладающем
$S$-свойством~\eqref{aspro2}.

Вторая теорема Шталя доказывается по следующей схеме.
Вначале на основе достаточно простых геометрических соображений
(ср. лемма~\ref{l1})
устанавливается, что существует допустимый компакт $F\in\KK_f$ такой, что
\begin{equation}
\mcap{(E,F)}=\min_{K\in\KK_f}\mcap{(E,K)}.
\label{asta2-3}
\end{equation}
Затем с помощью вариационного метода доказывается, что $F$ является
замыканием критических траекторий некоторого квадратичного дифференциала.
Отсюда уже вытекает $S$-свойство~\eqref{aspro2}. Непосредственно на основе
этого $S$-свойства устанавливается предельное соотношение~\eqref{ameco2} для
нулей знаменателей рациональных аппроксимаций $R_n$.
Соотношение~\eqref{ameco2} влечет сходимость по емкости последовательности
$\{R_n\}$ к функции~$f$ внутри (на компактных подмножествах) области $D$.
Единственность компакта $F$, удовлетворяющего условию~\eqref{mincon}, вытекает
из сходимости~\eqref{ameco2}.

Соотношение~\eqref{ameco2} эквивалентно соотношению
\begin{equation}
\bigl|Q_n(z)\bigr|^{1/n}\overset\mcap\longrightarrow
e^{-V^{\myt\lambda}(z)},\qquad z\in D,
\label{poco2}
\end{equation}
где старший коэффициент полинома $Q_n$ равен единице.


\end{document}